\documentclass[a4,10pt,fullpage]{article}
\usepackage{graphics}
\usepackage{graphicx}
\usepackage{epsfig}
\usepackage{amsmath}
\usepackage{amssymb}
\usepackage{fancybox}
\usepackage{amsfonts,url,times}
\usepackage{amscd}

\def\newmathop#1{\expandafter\gdef\csname #1\endcsname{\mathop{\rm #1}\nolimits}}
\def\newvmathop#1{\expandafter\gdef\csname v#1\endcsname{\mathop{\rm #1}\nolimits}}

 \makeindex
\newtheorem{corollary}{Corollary}
\newtheorem{definition}{Definition}
\newtheorem{theorem}{Theorem}
\newtheorem{proposition}{Proposition}
\newtheorem{rmk}{Remark}
\newtheorem{lemma}{Lemma}

\newvmathop{Hom}
\newvmathop{ord}
\newvmathop{Re}
\newvmathop{GL}
\newvmathop{Frob}
\newvmathop{PGL}
\newvmathop{Disc}
\newvmathop{Res}
\newvmathop{Tr}
\newvmathop{Norm}
\newvmathop{if}
\newvmathop{iff}
\newvmathop{and}
\newvmathop{or}
\newvmathop{for}
\newvmathop{sign}
\newvmathop{rk}
\newvmathop{corank}
\newvmathop{coker}
\newvmathop{codim}
\newvmathop{Ind}
\newvmathop{Res}
\newvmathop{tw}
\newmathop{Gal}
\newmathop{Aut}
\newmathop{Sel}
\newvmathop{Gal}
\newvmathop{Aut}
\newvmathop{Sel}

\newcommand{\qed}{~\rule{1.8mm}{1.8mm}}
\newcommand{\p}{\mathfrak{p}}
\newcommand{\q}{\mathfrak{q}}
\newcommand{\m}{\mathfrak{m}}

\newcommand{\isom}{\stackrel{\sim}{\rightarrow}}
\newcommand{\la}{\lambda}

\newcommand{\Z}{\mathbb{Z}}
\newcommand{\C}{\mathbb{C}}
\newcommand{\R}{\mathbb{R}}
\newcommand{\Q}{\mathbb{Q}}
\newcommand{\T}{\mathbb{T}}
\newcommand{\tH}{\mathbb{H}}
\newcommand{\ri}{\mathcal{O}} 
\newcommand{\A}{\mathbb{A}}
\newcommand{\E}{\textbf{E}}

\newcommand{\f}{\mathbf{f}}
\newcommand{\g}{\textbf{g}}

\newlength{\mpagewidth}
{\end{minipage}\end{Sbox}\fbox{\TheSbox}\medbreak}

\begin{document}

\title{Special values of L-functions and false Tate curve extensions II}
\author{Thanasis Bouganis}
\maketitle
\begin{abstract}
In this paper we show how one can combine the $p$-adic
Rankin-Selberg product construction of Hida with freeness results of
Hecke modules of Wiles to establish interesting congruences between
special values of L-functions. These congruences is a part of some
deep conjectural congruences that follow from the work of Kato on
the non-commutative Iwasawa theory of the false Tate curve
extension.
\end{abstract}

\section{Introduction}

Let $E$ be an elliptic curve defined over $\Q$ and $p$ a rational
prime. In the classical setting of cyclotomic Iwasawa theory for
elliptic curves one is concerned with the study of the twists of the
elliptic curve by finite order character that factor through the
cyclotomic $\Z_{p}$ extension $\Q_{cyc} \subset \cup_{n \geq
0}\Q(\mu_{p^{n}})$, where $\mu_{p^{n}}$ is the group of the
$p^{n}$-th roots of unity. The aim of the theory is to obtain a link
between the analytically defined $L$ functions attached to $E$, and
its twists, and the arithmetic properties of the elliptic curve over
the cyclotomic tower. The cyclotomic Main Conjecture for elliptic
curves gives to this conjectural link a very precise form. We note
that much has already been proven towards this Main Conjecture by
Kato \cite{Kato2}, and Skinner and Urban have announced a complete
proof for semi-stable $E$, subject to proving certain results about
the Galois representations attached to automorphic forms.

One of the key ingredients of the above Main Conjecture are the
$p$-adic $L$ functions. These are usually realized as $p$-adic
measures over Galois groups, which, when evaluated at finite order
characters, interpolate canonically modified values of the $L$
function. Their construction usually involves two steps. The first
one is to find proper transcendental numbers, usually called
periods, such that the ratio of the $L$ values over these periods
gives an algebraic number. The second step is to prove that these
values, or a slight modification of them, have the desired
interpolation and integrality properties.

Lately there has been great interest in extending the classical
Iwasawa theory to a non abelian setting, that is to replace the
$\Z_{p}$ extension by more general $p$-adic Lie extensions whose
Galois group is non-abelian. In fact in \cite{CFKSV} a precise
analogue of the Main Conjecture in this non abelian setting for a
large family of $p$-adic Lie groups has been stated.

One of the extensions that is of particular interest is the so
called ``false Tate curve'' extensions. That is extensions of the
form, $\Q_{FT}:= \cup_{n \geq 0} \Q(\mu_{p^{n}},\sqrt[p^{n}]{m})$
for some $p$-power free integer $m>1$. Note that the Galois group is
the semi-direct product $\Z_{p} \ltimes \Z_{p}^{\times}$. There is a
conjectural theory for $p$-adic $L$ functions that should exist in
this setting. In a work with V.Dokchitser \cite{BougDok} we have
addressed the first of the above mentioned two steps, that is
algebraicity of the critical values of the $L$ functions involved.

In order to make things more explicit let us fix some more notation.
We write $E$ for an elliptic curve defined over $\Q$ and $N_{E}$ for
its conductor. As we already mentioned we consider the extensions
$\Q_{FT,n} := \Q(\mu_{p^{n}},\sqrt[p^{n}]{m})$ and $\Q_{FT}=\cup_{n
\geq 0}\Q_{FT,n}$. We write $\rho$ for an Artin representation that
factors through $\Q_{FT}$ and $N_{\rho}$ for its conductor. Let us
also write $L(E,\rho,s)$ for the $L$ function attached to $E$
twisted by $\rho$. We consider the value of $L(E,\rho,s)$ at the
critical point $s=1$. The fact that the Artin representations $\rho$
factor through the false Tate curve allowed us to establish the
analyticity of $L(E,\rho,s)$ at $s=1$ and then our main result in
\cite{BougDok} is concerned with the algebraic properties of these
values. Let us write $\Omega_{\pm}(E)$ for the N\'{e}ron periods
attached to the elliptic curve $E$. Then we have shown that
\[
\frac{L(E,\rho,1)}{\Omega_{+}(E)^{dim(\rho^{+})}\Omega_{-}(E)^{dim(\rho^{-})}}
\in \overline{\Q}.\]

\noindent for all Artin representations $\rho$ that factor through
$\Q_{FT}$. Actually we did more. Namely, involving also the period
that should correspond to the ``Artin motive" $M(\rho)$ attached to
$\rho$ we established the period conjecture of Deligne that gives a
precise description of the number field where this value lies.

Let us now move to the second step that we mentioned above, that is
the $p$-adic properties of these values. From now on we will assume
that the elliptic curve has good ordinary reduction at $p$. We start
by stating a conjectural congruence between these $L$ values for
different Artin representations. We define the quantity $R(\rho)$ as
\[
R(\rho) :=
e_{p}(\rho)u^{-v_{p}(N_{\rho})}\frac{P_{p}(\hat{\rho},u^{-1})}{P_{p}(\rho,w^{-1})}
\cdot \frac{L_{\{p,q|m\}}(E, \rho,1)}{\Omega_{+}(E)^{dim
(\rho^{+})}\Omega_{-}(E)^{dim (\rho^{-})}}
\]
\noindent where $e_{p}(\rho)$ is a local epsilon factor of $\rho$
suitably normalized, $P_{p}(\rho,X)$ is the usual characteristic
polynomial associated to $\rho$ at $p$ and $u,w$ are $p$-adic
numbers defined by,
\[
1-a_{p}X+pX^{2} = (1-uX)(1-wX),\,\,\,\, u\in
\Z_{p}^{\times}\,\,\,and\,\,\,p+1-a_{p}=\#E_{p}(\mathbb{F}_{p})
\]
\noindent Here $\hat{\rho}$ is the dual representation but in our
false Tate curve setting it is easy to see that $\hat{\rho}\cong
\rho$. Finally the subscript $\{p,q|m\}$ means that we have removed
the Euler factors at these primes. Then we state,

\paragraph{Conjecture:}For each $n \geq 1$, let $\chi_{n}$ be a character
of $Gal(\Q_{FT,n}/\Q(\mu_{p^{n}}))$ of exact order $p^{n}$. Write
$\rho_{n}$ for the induced representation of $\chi_{n}$ to
$Gal(\Q_{FT,n}/\Q)$ and $\sigma_{n}$ for the representation induced
to $Gal(\Q_{FT,n}/\Q)$ from the trivial one over $\Q(\mu_{p^{n}})$.
Then, the values $R(\rho_{n})$ and $R(\sigma_{n})$ are $p$-adically
integral and satisfy
\[
|R(\rho_{n}) - R(\sigma_{n})|_{p} < 1
\]
\noindent or more generally
\[
|R(\rho_{n} \otimes \psi) - R(\sigma_{n} \otimes \psi)|_{p} < 1
\]
\noindent where $\psi$ is a finite order character of
$Gal(\Q^{cyc}/\Q)$ and $|\cdot|_{p}$ normalized as $|p|_{p} =
p^{-1}$.

Let us comment a little bit more on this conjecture and its
connection to non commutative Iwasawa theory. The definition of the
quantity $R(\rho)$ describes the interpolation properties that the
conjectural, as in \cite{CFKSV}, non-abelian $p$-adic $L$-function
should satisfy. Indeed the authors in \cite{CFKSV} have conjectured
the existence of an element in the $K_{1}$ of the Iwasawa algebra
associated to this extension that interpolates suitably modified, as
above, values of $L(E,\rho,1)$ and plays the role of the non-abelian
$p$-adic $L$ function in their theory. Note that the representations
$\rho_{n}$ and $\sigma_{n}$ are defined over $\Q$ and are congruent
modulo $p$ that is if we consider their reduction modulo $p$ then
their semi-simplifications are isomorphic. Hence the existence of
the non-abelian $p$-adic $L$ function would imply that its values
should be also $p$-adically close.

There is almost nothing known concerning the construction of this
object for a general $p$-adic Lie extension. However in the setting
that we are interested in, the false Tate curve extension, Kato in
\cite{Kato1} has related the existence of this non-abelian object
with congruences between classical abelian $p$-adic $L$ functions
over various fields of the extension. We take some time to explain
this as it will help us motivate the results that appear in this
paper. Let $G$ be the Galois group of the false Tate curve extension
and $\Lambda(G)=\Z_{p}[[G]]$ the Iwasawa algebra of $G$. We set
$U^{(n)}:= ker (\Z_p^{\times}\rightarrow (\Z/p^n\Z)^{\times})$. The
main result of Kato in \cite{Kato1} is the construction of an
injective homomorphism
\[
\theta_G : K_1(\Lambda(G)) \rightarrow \prod_{n \geq
0}\Z_p[[U^{(n)}]]^{\times}
\]
\noindent and the explicit description of the image. In order to
make this last statement a little bit more precise we write, for $n
\geq m \geq 0$, $N_{m,n}:\Z_p[[U^{(m)}]]\rightarrow \Z_p[[U^{(n)}]]$
for the canonical norm map, $\phi$ be the ring homomorphism
$\Z_p[[\Z_p^{\times}]]\rightarrow \Z_p[[\Z_p^{\times}]]$ induced by
the rising to the power $p$ map on $\Z_p^{\times}$. Then the result
of Kato says that $\theta_G(K_1(\Lambda(G)))=(a_n)_{n \geq 0}$ with
\[
\prod_{0 < i \leq n}N_{i,n}(c_i)^{p^{i}} \equiv 1 \mod{p^{2n}}
\]
\noindent with $c_n = b_n \phi(b_{n-1})^{-1}$ and $b_n=a_n
N_{0,n}(a_0)^{-1}$. The elements $a_n$ have an arithmetic meaning,
they are abelian $p$-adic $L$ functions. More precisely if we write
$\rho_n$ for the Artin representation of $G$ induced from a
character of $p^{n}$ order of the Galois group
$Gal(\Q(\mu_{p^{n}},\sqrt[p^{n}]{m})/\Q(\mu_{p^n}))$, then the
elements $a_n$ are the abelian $p$-adic $L$-functions interpolating
the values $L(E\otimes \rho_n \otimes \chi,1)$, for $\chi$ Dirichlet
characters of the cyclotomic extension of $\Q$.

The conjectural congruences that we have written above correspond to
the case of $n=1$ of Kato's congruences after evaluating the abelian
$p$-adic $L$ functions at the character $\psi$. There is
computational support for these conjectures; initially by Balister
\cite{Balister} and much more vastly by the Dokchitser brothers
\cite{TDVD}. In the first part of this work \cite{Boug} we have
showed the existence of the abelian $p$-adic $L$-functions $a_{n}$
appeared in Katos's congruences and proved the above conjectural
congruences up to an issue of periods. Namely there we have used not
the motivic periods that are stated in the congruences but
automorphic periods, the so called Eichler-Shimura-Harder periods,
that appear quite natural in the so called modular symbol
construction. There we came across to a rather deep problem, namely
the relation of these automorphic periods as one use the functorial
properties of the $L$-functions and especially base-change. We say a
little bit more on this at the last section of this paper. Finally
we note that in \cite{DelbWard} an inductive argument was used to
show how these congruences $($for $n=1)$ can provide congruences for
$n>1$ in the form conjectured by Kato but unfortunately not modulo
the right $p$ power.

Our aim in this paper is to tackle the conjectural congruences
insisting on getting the right motivic periods. We achieve that for
the case where $p=3$ but we also discuss possible extensions for the
case of $p>3$. We need to impose some further conditions on $E$,
other of technical nature which we believe can be removed and other
that seem important. Namely from now on we assume that (a) The curve
$E$ is semi-stable and if we consider the minimal discriminant
$\Delta_{E}=\prod_{q|N_{E}}q^{i_{q}}$ then $p$ does not divide
$i_{q}$ for all $q$. Note that the last condition means that the
conductor of $E$ is equal to the Artin conductor of the mod $p$
representation obtained by $E$. (b) We assume that $m$ that appear
in the false Tate extension is power free with $(m,N_{E})=(m,p)=1$
and, (c) a rather important assumption, that $E$ has no rational
subgroup of order $p$, that is the associated modulo $p$
representation is irreducible. Finally we mention here that as our
aim here is to address the issue of motivic versus automorphic
periods we focus on proving the above conjectures for $\psi=1$.
However we lay all important constructions so that everything can be
extend to the case $\psi$ being not trivial.

Our proof can be divided into two parts. Let us write $f \in
S_{2}(\Gamma_{0}(N_{E});\Q)$ for the rational newform that we can
associate to $E$. In the first part we rely on the work of Hida of
the construction of a p-adic Rankin-Selberg product initiated in
\cite{Hida2} and generalized in \cite{Hida4}. We can associate a
newform $g$ of weight one to the Artin representation $\rho$ and an
Eisenstein series $\mathcal{E}$ of weight one with $\sigma$. Using
them, we construct $p$-adic measures $d\mu_{f,g}$ and
$d\mu_{f,\mathcal{E}}$ over $\Z_{p}^{\times}$ that are congruent
modulo $p$, in the sense that their values at every finite character
of $\Z_{p}^{\times}$ are congruent. These measures interpolate,
$p$-adically, twists of the critical values of the Rankin-Selberg
products $D(f,g,s)$ and $D(f,\mathcal{E},s)$ by finite order
characters. Evaluating these measures at the trivial character we
get a first form of congruences between $D(f,g,1)$ and
$D(f,\mathcal{E},1)$. Under the semi-stable assumption we can easily
relate the Rankin-Selberg product to the twists of the elliptic
curve $E$.

However we do not yet get the congruences stated in the theorem
above. We need to work further two things. First, in order to
establish the congruences between the measures above, we had to
clear a denominator $c(f,m)$ that depends solely on $f$ and $m$.
Hence we get congruences after multiplying with this constant
$c(f,m)$. Second, the periods that we use to get the rationality of
the Rankin-Selberg product are closely related to the Petersson
inner product $<f,f>$. These periods may not be equal to our periods
$\Omega_{+}(E)$ and $\Omega_{-}(E)$ up to a $p$-adic unit. These two
problems are related. That is, the reason that the denominator
$c(f,m)$ appears in our $p$-adic interpolation is the fact that the
Petersson inner product is not the proper automorphic period in
order to get $p$-adically integral ratios of the form
$\frac{L-values}{aut.\,\,periods}$.

In the second part we show, under the assumptions of the theorem,
that indeed this is the case. This part relies heavily on the work
of Wiles. We make use of two of his important results in
\cite{Wiles}. The first one is an extension of a theorem of Mazur
\cite{Mazur1} on the freeness, over a completed Hecke algebra, of
the first cohomology group of modular curves after localizing it at
a proper maximal ideal. The second one is an extension of a theorem
of Ihara on the study of maps between Jacobians of modular curves of
different levels. Here we would like to mention how helpful was for
us the paper of Darmon, Diamond and Taylor \cite{DDT} reviewing the
work of Wiles.

Let us just mention that we tried to apply the same ideas for $p>3$.
Here in order to bring things to the previous setting we use the
fact that the base-change property for automorphic representations
of $GL(2)$ has been proved for cyclic extensions \cite{Langlands}.
Using this, we can work the congruences over the totaly real field
$F := \Q(\mu_{p})^{+}$. However we face two problems. First the fact
that we work with a prime that ramifies in $F$ puts restrictions on
the freeness results that we need. Second we need to relate our
defined automorhic periods over $F$ with the ones over $\Q$, and
even stronger we need the relation to be up to $p$-adic units a
problem much of the same nature that we face in our work
\cite{Boug}. We do not have an answer to these questions yet.

\paragraph{Acknowledgements:} The author would like to thank
Professor John Coates for suggesting to work on Kato's congruences
and for recommending to consider the use of the Rankin-Selberg
method and its $p$-adic version.

\section{Basic Notations}
Let $\tH$ be the complex upper half plane. If we denote by
$GL_{2}^{+}(\R)$ the two by two real matrices with positive
determinant, then we consider the action of them on $\tH$ by liner
fractional transformations, $z \mapsto \alpha(z)=\frac{az+b}{cz+d}$,
for
$\alpha = \left(%
\begin{array}{cc}
  a & c    \\
  b & d \\
\end{array}%
\right) \in GL_{2}^{+}(\R)$. We let $k \geq 1$ be an integer and we
define an action of $GL_{2}^{+}(\R)$ on functions $f:\tH \rightarrow
\C$ by
\[
f \mapsto (f|_{k}[\alpha])(z) = det(\alpha)^{k/2}
(cz+d)^{-k}f(\alpha(z))
\]
for $\alpha = \left(%
\begin{array}{cc}
  a & c    \\
  b & d \\
\end{array}%
\right) \in GL_{2}^{+}(\R)$. We denote by $SL_{2}(\Z)$ the two by
two matrices with determinant 1 and integral entries. For a positive
integer $N$ we have the standard notations for the subgroups of
$SL_{2}(\Z)$,
\[
\Gamma(N) = \{\gamma \in SL_{2}(\Z) \,\,\,|\,\,\, \gamma \equiv \left(%
\begin{array}{cc}
  1 & 0    \\
  0 & 1 \\
\end{array}%
\right) \mod{N} \}
\]
\[
\Gamma_{0}(N) = \{\gamma \in SL_{2}(\Z) \,\,\,|\,\,\, \gamma \equiv \left(%
\begin{array}{cc}
  * & *    \\
  0 & * \\
\end{array}%
\right) \mod{N} \}
\]
\[ \Gamma_{1}(N) = \{\gamma \in \Gamma_{0}(N) \,\,\,|\,\,\, \gamma \equiv \left(%
\begin{array}{cc}
  1 & *    \\
  0 & 1 \\
\end{array}%
\right) \mod{N} \} \]

We write $M_{k}(\Gamma_{1}(N))$ (resp. $S_{k}(\Gamma_{1}(N))$) for
the space of modular forms (resp. cusp forms) of weight $k$ with
respect to $\Gamma_{1}(N)$. We write $M_{k}(\Gamma_{0}(N),\chi)$
(resp $S_{k}(\Gamma_{0}(N),\chi)$ for modular forms (resp. cusp
forms) with respect to $\Gamma_{0}(N)$ and Nebentype $\chi$.

Let us consider a cusp form $f \in S_{k}(\Gamma_{0}(N),\chi)$ and a
modular form $g \in M_{l}(\Gamma_{0}(N),\psi)$, for some integers
$k$ and $l$ where we moreover assume $k>l$. Let us write their
Fourier expansions at $\infty$ cusp as $f(z)=\sum_{n=1}^{\infty}
a(n,f)q^{n}$ and $g(z)=\sum_{n=0}^{\infty}a(n,g)q^{n}$ with
$q=e^{2\pi \imath z}$. We also define
$f^{\rho}(z)=\sum_{n=1}^{\infty}\overline{a(n,f)}q^{n} \in
S_{k}(\Gamma_{0}(N),\bar{\chi})$. We consider the quantities
$L(f,g,s):=\sum_{n=1}^{\infty}a(n,f)a(n,g)n^{-s}$ and their
Rankin-Selberg convolution, $D(f,g,s):=
L_{N}(\chi\psi,2s+2-k-l)L(f,g,s)$ where we have removed the Euler
factors at $N$ from $L(\chi\psi,s)$. If we assume that $f$ and $g$
are actually normalized eigenforms and if we write their $L$
functions
$L(f,s)=\prod_{q}\{(1-\alpha(q,f)q^{-s})(1-\beta(q,f)q^{-s})\}^{-1}
$ and
$L(g,s)=\prod_{q}\{(1-\alpha(q,g)q^{-s})(1-\beta(q,g)q^{-s})\}^{-1}$
then we have that
\[
D(f,g,s)=
\prod_{q}\{(1-\alpha(q,f)\alpha(q,g)q^{-s})(1-\alpha(q,f)\beta(q,g)q^{-s})\times\]
\[(1-\beta(q,f)\alpha(q,g)q^{-s})(1-\beta(q,f)\beta(q,g)q^{-s})\}^{-1}
\]

\section{$p$-adic modular forms and measures}
In this section we introduce the needed background in order to
obtain the $p$-adic version of the Rankin-Selberg convolution. For
all this background we follow Hida's papers \cite{Hida2,Hida4}. We
let $p$ be a prime number and we fix an embedding $\overline{\Q}
\hookrightarrow \overline{\Q}_{p} \hookrightarrow \C_{p}$, where
$\C_{p}$ is the $p$-adic completion of $\overline{\Q}_{p}$ under the
normalized $p$-adic absolute value $|\cdot|_{p}$ with $|p|_{p} =
p^{-1}$. For any subring $R \subseteq \overline{\Q}$ we consider the
$R$-modules,
\[
M_{k}(\Gamma_{0}(N),\psi;R) := \{f \in M_{k}(\Gamma_{0}(N),\psi)
\mid \,\,\,\,\, f(z) = \sum_{n \geq 0} a(n,f)q^{n}, \,\,\,a(n,f) \in
R \}
\]
\[
M_{k}(\Gamma_{1}(N);R) := \{f \in M_{k}(\Gamma_{1}(N)) \mid
\,\,\,\,\,\, f(z) = \sum_{n \geq 0} a(n,f)q^{n},\,\,\,a(n,f) \in R
\}
\]
Moreover we define $S_{k}(\Gamma_{0}(N),\psi;R) =
S_{k}(\Gamma_{0}(N),\psi) \cap M_{k}(\Gamma_{0}(N),\psi;R)$ and
similar for $S_{k}(\Gamma_{1}(N);R)$. For a modular form $f \in
M_{k}(\Gamma_{1}(N);\overline{\Q})$ it is known that one can define
the $p$-adic norm of $f$, $|f|_{p} := sup_{n \geq 0} |a(n,f)|_{p}$.
Let now $K_{0}$ be any finite extension of $\Q$ and write $K$ for
the closure of $K_{0}$ in $\C_{p}$. We define the space
$M_{k}(\Gamma_{0}(N),\psi;K)$ (resp. $M_{k}(\Gamma_{1}(N);K)$) to be
the $p$-adic completion of the space
$M_{k}(\Gamma_{0}(N),\psi;K_{0})$ (resp.
$M_{k}(\Gamma_{1}(N);K_{0})$ with respect to the norm $|\cdot|_{p}$
inside $K[[q]]$ where we consider $q$ as indeterminant. Then it is
known by the work of Deligne and Rapoport \cite{DeligneRapoport}
that, $M_{k}(\Gamma_{0}(N),\psi;K) = M_{k}(\Gamma_{0}(N),\psi;K_{0})
\otimes_{K_{0}} K$, $M_{k}(\Gamma_{1}(N);K) =
M_{k}(\Gamma_{1}(N);K_{0}) \otimes_{K_{0}} K$. Moreover it is known
that the definition of $M_{k}(\Gamma_{1}(N);K)$ and
$M_{k}(\Gamma_{0}(N),\psi;K)$ is independent of the choice of the
dense subfield $K_{0}$. Let us now write $\ri_{K}$ for the $p$-adic
ring of integers of $K$. Then we define the $p$-adic integral
modular forms as,
\[
M_{k}(\Gamma_{0}(N),\psi;\ri_{K}):=\{f \in
M_{k}(\Gamma_{0}(N),\psi;K) \mid |f|_{p} \leq 1 \} =
M_{k}(\Gamma_{0}(N),\psi;K) \cap \ri_{K}[[q]],
\]
\[
M_{k}(\Gamma_{1}(N);\ri_{K}):=\{f \in M_{k}(\Gamma_{1}(N);K) \mid
|f|_{p} \leq 1 \} = M_{k}(\Gamma_{1}(N);K) \cap \ri_{K}[[q]]
\]
\begin{definition} $(p$-adic modular forms$)$. Let $A$ be either $K$
or $\ri_{K}$. We consider the spaces,
\[
M_{k}(N;A) := \cup_{n=0}^{\infty} M_{k}(\Gamma_{1}(Np^{n});A)
\,\,\,and \,\,\, M_{k}(N,\psi;A) := \cup_{n=0}^{\infty}
M_{k}(\Gamma_{0}(Np^{n}),\psi;A)
\]
Then we define the space of $p$-adic modular forms of
$\Gamma_{1}(N)$, resp. of $\Gamma_{0}(N)$ and character $\psi$, as
the completion of the above spaces with respect to the norm
$|\cdot|_{p}$. We denote them by $\overline{M}_{k}(N;A)$, resp.
$\overline{M}_{k}(N,\psi;A)$.
\end{definition}

We note that all the above discussion can be done considering cusp
forms instead of modular forms. In particular we can consider also
$p$-adic cusp forms which we will denote by $\overline{S}_{k}(N,A)$
and $\overline{S}_{k}(N,\psi;A)$.

\begin{rmk} \label{rmk:weight}
For our later use, we mention that the space $\overline{M}_{k}(N,A)$
is actually independent of $k$ for $k \geq 2$, so we may also write
just $\overline{M}(N;A)$, see \cite{Hida4}.
\end{rmk}

Now we are going to define $p$-adic Hecke operator that extend the
usual ones when restricted to the space of classical modular forms.
For any integer $n$ prime to $N$ we consider a matrix $\sigma_{n}
\in
\Gamma_{0}(N)$, such that $\sigma_{n} \equiv \left(%
\begin{array}{cc}
  n^{-1} & o \\
  0 & n \\
\end{array}%
\right) \mod{N}$. It follows by the work of Deligne and Rapoport
\cite{DeligneRapoport} that the action $f \mapsto f |_{k}
\sigma_{n}$ on $M_{k}(\Gamma_{1}(N);K)$ is integral, that is it
preserves the integral space $M_{k}(\Gamma_{1}(N);\ri_{K})$. We
``define'' the Hecke operators $T(\ell)$ and $S(\ell)$, for every
prime $\ell$, acting on $M_{k}(\Gamma_{1}(N);K)$ by describing their
action on the $q$-expansion,
\[
a(n,T(\ell)f) = \left\{%
\begin{array}{ll}
    a(\ell n,f)+\ell^{k-1}a(\frac{n}{\ell},f|_{k}\sigma_{\ell}), & \hbox{if $\ell$ is prime to $N$;} \\
    a(\ell n,f), & \hbox{otherwise.} \\
\end{array}%
\right.
\]
\[
a(n,S(\ell)(f)) = \left\{%
\begin{array}{ll}
    \ell^{k-2}a(n,f|_{k}\sigma_{\ell}), & \hbox{if $\ell$ is prime to $N$;} \\
    0, & \hbox{otherwise.} \\
\end{array}%
\right.
\]
Note that these definitions are consistent with the ones on the
classical elliptic modular forms. We define the Hecke algebra
$H_{k}(\Gamma_{0}(N),\psi;A)$, resp. $H_{k}(\Gamma_{1}(N);A))$, for
$A$ either $K$ or $\ri_{K}$ as the $A$-subalgebra of
$End_{A}(M_{k}(\Gamma_{0}(N),\psi;A))$, resp.
$End_{A}(M_{k}(\Gamma_{1}(N);A))$, generated by $T(\ell)$ and
$S(\ell)$ for all primes $\ell$. Similarly we define
$h_{k}(\Gamma_{0}(N);\psi;A)$ and $h_{k}(\Gamma_{1}(N);A)$ when we
restrict the action to the space of cusp forms. Actually one has
that $H_{k}(\Gamma_{0}(N),\psi;A)=H_{k}(\Gamma_{0}(N),\psi;\Z)
\otimes_{\Z} A$ and similarly for the other spaces. Finally we note
that when $p|N$ the action of the Hecke operators is $p$-adically
integral i.e. $|Tf|_{p} \leq |f|_{p}$ for every $T \in
H_{k}(\Gamma_{1}(N);\ri_{K})$.

We now define $p$-adic Hecke algebras. Notice that we have the
$\ri_{K}$-surjective homomorphisms induced by restriction of the
Hecke operators,
\[
H_{k}(\Gamma_{0}(Np^{m}),\psi;\ri_{K}) \rightarrow
H_{k}(\Gamma_{0}(Np^{n}),\psi;\ri_{K}) \,\,\,\,\,for \,\,\,m \geq n
\geq 1
\]
\[
H_{k}(\Gamma_{1}(Np^{m});\ri_{K}) \rightarrow
H_{k}(\Gamma_{1}(Np^{n});\ri_{K}) \,\,\,\,\,for \,\,\,m \geq n \geq
1
\]
\begin{definition}We define the space of $p$-adic Hecke algebras
$H_{k}(N,\psi;\ri_{K})$ (resp. $H_{k}(N;\ri_{K}))$ by the projective
limit, $\underleftarrow{\lim}_{n}H_{k}(\Gamma_{0}(N),\psi;\ri_{K})$
(resp. $\underleftarrow{\lim}_{n}H_{k}(\Gamma_{1}(N);\ri_{K}))$.
Similarly we define the spaces $h_{k}(N,\psi;\ri_{K})$ and
$h_{k}(N;\ri_{K})$.
\end{definition}

By definition this operators act on the spaces $M_{k}(N;A)$ and
$M_{k}(N,\psi;A)$ for $A$ equal to $K$ or $\ri_{K}$. However the
fact they are $p$-adically integral allow us to extend their action
to the space of $p$-adic modular forms $\overline{M}_{k}(N;A)$ and
$\overline{M}_{k}(N,\psi;A)$. Our next step is to define Hida's
ordinary idempotent $e$ attached to the Hecke operator $T(p)$. We
start with a general lemma,

\begin{lemma} For any commutative $\ri_{K}$-algebra $R$ of finite rank over $\ri_{K}$ and for any $x
\in R$ the limit $\lim_{n \rightarrow \infty}x^{n!}$ exists and
gives an idempotent of $R$.
\end{lemma}

\paragraph{Proof} See \cite{Hida7} (p.201) \qed.

\begin{definition} We define an idempotent $e_{n}$ in
$H_{k}(\Gamma_{0}(Np^{n},\psi;\ri_{K})$ and in
$H_{k}(\Gamma_{1}(Np^{n};\ri_{K})$ by the limit $e_{n} = \lim_{m
\rightarrow \infty} T(p)^{m!}$. Moreover we define an idempotent in
$H_{k}(N;\ri_{K})$ and in $H_{k}(N,\psi;\ri_{K})$ by taking the
projective limit $e = \underleftarrow{\lim}_{n} e_{n}$.
\end{definition}

We will be interested in the space
$e\overline{M}_{k}(N,\psi;\ri_{K})$, usually called the ordinary
part of $\overline{M}_{k}(N,\psi;\ri_{K})$ and denoted by
$\overline{M}_{k}^{\circ}(N,\psi;\ri_{K})$. Actually this space is
not that large as the following lemma indicates,

\begin{lemma}\label{lemma:projector}$($Hida$)$ Let $C(\psi)$ be the conductor of the character
$\psi$. Define positive integers $N'$ and $C(\psi)'$ by writing
$N=N'p^{r}$ and $C(\psi)=C(\psi)'p^{t}$ with
$(N',p)=(C(\psi)',p)=1$. Let $s:= max(t,1)$. Then,
\[
e\overline{M}_{k}(N,\psi;\ri_{K}) \subset
M_{k}(\Gamma_{0}(N'p^{s}),\psi;\ri_{K})
\]
\end{lemma}
\paragraph{Proof:} See \cite{Hida2}.

\begin{definition}
We say that a normalized eigenform $f_{0} \in
S_{k}(\Gamma_{0}(N_{0})\psi)$ is an ($p$-) ordinary form if,

\begin{enumerate}
    \item The level $N_{0}$ of the form $f$ is divisible by $p$.
    \item The Fourier coefficient $a(p,f_{0})$ is a $p$-adic unit.
\end{enumerate}
\end{definition}

The following lemma is proved in \cite{Hida2} (p. 168),

\begin{lemma} \label{lemma:ordinary}Let $f \in S_{k}(\Gamma_{0}(N),\psi)$ be a newform
with $k \geq 2$ and $|a(p,f)|_{p}=1$. Then, there is a unique
ordinary form $f_{0}$ of weight $k$ and character $\psi$ such that
$a(n,f) = a(n,f_{0})$ for all $n$ not divisible by $p$. Moreover,
$f_{0}$ is given by,
\[
f_{0}(z) = \left\{%
\begin{array}{ll}
    f(z), & \hbox{if $p$ divides $N$;} \\
    f(z)-w f(pz), & \hbox{otherwise.} \\
\end{array}%
\right.
\]
where $w$ is the unique root of $X^{2}-a(p,f)X+\psi(p)p^{k-1}=0$
with $|w|_{p}<1$. Moreover in the second case i.e. $(p,N)=1$ we have
that $N_{0}=Np$ and that $a(p,f_{0}) = u$ where $u$ is the $p$-adic
unit root of the above equation.
\end{lemma}

Let us now consider a surjective $K$-linear homomorphism
$\Phi:h_{k}(\Gamma_{0}(N_{0}),\psi;K) \rightarrow K$ that is induced
by an ordinary form $f_{0}$ by sending $T(n) \mapsto a(n,f_{0})$.
Let us moreover assume that this map is split (we will show later
that in the case of interest this will be true) and induces an
algebra direct decomposition, $h_{k}(\Gamma_{0}(N_{0}),\psi;K) \cong
K \times A$ for some summand $A$ and let us denote by $1_{f_{0}}$
the idempotent corresponding to the first summand isomorphic to $K$.
We now consider the linear form
$\ell_{f_{0}}:\overline{S}_{k}(N_{0},\psi;K) \rightarrow K$ defined
by, $\ell_{f_{0}}(g) := a(1, 1_{f_{0}}e\,\, g)$. Note that, by lemma
\ref{lemma:projector}, the linear form is well defined.

\begin{proposition}$($Hida's linear operator$)$\label{proposition:formula} Assume that $K_{0}$ contains all the Fourier
coefficients of the ordinary form $f_{0}$. Then, the linear form
$\ell_{f_{0}}$ has values in $K_{0}$ on
$S_{k}(\Gamma_{0}(N_{0}p^{n}),\psi;K_{0})$ for every $n \geq 0$.
Furthermore, for $g \in S_{k}(\Gamma_{0}(N_{0}p^{n}),\psi;K_{0})$ we
have
\[
\ell_{f_{0}}(g)=a(p,f_{0})^{-n}p^{n(k/2)}\frac{<h_{n},g>_{N_{0}p^{n}}}{<h,f_{0}>_{N_{0}}}
\]
where $h=f_{0}^{\rho}|_{k}\left(%
\begin{array}{cc}
  0 & -1    \\
  N_{0} & 0 \\
\end{array}%
\right)$, $h_{n}(z)=h(p^{n}z)$.
\end{proposition}

\paragraph{Proof}See \cite{Hida2} p.175.

\paragraph{}
We note that if we consider a constant $c(f_{0}) \in \ri_{K}$ such
that $c(f_{0})1_{f_{0}} \in h_{k}(\Gamma_{1}(N_{0});\ri_{K})$ then
we have an integral valued linear form $c(f_{0})\ell_{f_{0}} :
\overline{S}_{k}(N_{0},\psi;\ri_{K}) \rightarrow \ri_{K}$ as the
Hecke operators are $p$-adically integral

\paragraph{$p$-adic modular forms valued measures:} Now we are going to define $p$-adic measures associated with
$p$-adic modular forms $\overline{M}(N;\ri_{K})$ for some $N$
relative prime to $p$. Note that it follows from remark 1 that we do
not need to specify the weight.

We let $X$ to be a $p$-adic space that consists of some copies of
$\Z_{p}$ and of a finite product of finite groups. For our
applications later $X$ is going to be just $\Z_{p}^{\times}\cong
(1+p\Z_{p}) \times (\Z/p\Z)^{\times}$. Let us write $C(X;\ri_{K})$
for the space of continuous functions of $X$ with values in
$\ri_{K}$ and $LC(X;\ri_{K})$ for the space of locally constant
functions on $X$. A measure $\mu$ on $X$ with values in the space
$\overline{M}(N;\ri_{K})$ is just an $\ri_{K}$-linear homomorphism
from $C(X;\ri_{K})$ to $\overline{M}(N;\ri_{K})$.

Let us consider the space $Z_{N}:=\Z_{p}^{\times} \times
(\Z/N\Z)^{\times}$ and for an element $z \in Z_{N}$ let us write
$z_{p}$ for the projection of $z$ to the first component. We can
define an action of $Z_{N}$ on the space
$M_{k}(\Gamma_{1}(Np^{r});\ri_{K})$ by $f \mapsto
f|z:=z_{p}^{k}f|_{k}\sigma_{z}$ with $\sigma_{z}$ as defined above.
This action can be extended to $\overline{M}(N,\ri_{K})$ (see
\cite{Hida4} p. 10).

\begin{definition} (see \cite{Hida4})
We say that a $p$-adic measure $\mu:C(X;O_{K}) \rightarrow
\overline{M}(N;O_{K})$ is arithmetic if the following three
conditions are satisfied,

\begin{enumerate}
    \item There exists positive integer $k$ such that for every $\phi
    \in LC(X;\ri_{K})$,
\[
\mu(\phi) \in M_{k}(Np^{\infty};\ri_{K})
\]
We will call $k$ the weight of $\mu$.
    \item There are continuous action $Z_{N} \times X \rightarrow
    X$ and a finite order character $\xi : Z_{N} \rightarrow
    \ri^{\times}_{K}$ such that
    $\mu(\phi)|z=z^{k}_{p}\xi(z)\mu(\phi(z\cdot x))$ for every $\phi
    \in C(X;O_{K})$, where $k$ the weight of $\mu$.
    We then say that the arithmetic measure is of character $\xi$.
\end{enumerate}

We say that the measure is cuspidal if $\mu$ actually takes values
in $\overline{S}(N;\ri_{K})$.
\end{definition}

We are interested in attaching arithmetic measures to a given
modular form. Given a modular form $f \in
M_{k}(\Gamma_{0}(N),\chi;\ri_{K})$ with $q$-expansion $f(z)=\sum_{n
\geq 0} a(n,f)q^{n}$ we can associate a measure $d\mu_{f}$ on
$X:=\Z_{p}^{\times}$ by,
\[
d\mu_{f}(\phi) \mapsto \sum_{n \geq 1} \phi(n)a(n,f)q^{n},
\,\,\,\,\,\phi \in C(X;\ri_{K})
\]
where we define the action of $Z_{N}$ on $\Z_{p}^{\times}$ by $z
\cdot x \mapsto z_{p}^{2}x$. From the following lemma due to Shimura
we conclude that $d\mu_{f}$ is an arithmetic measure of weight $k$
and character $\chi$.

\begin{lemma} Let $g=\sum_{n=0}^{\infty}b(n,g)q^{n} \in
M_{k}(\Gamma_{0}(N),\omega)$ and $\phi$ an arbitrary function on
$Y_{m}=$\textbf{Z}$/Np^{m}$\textbf{Z}. Define
$g(\phi):=\sum_{n=0}^{\infty}\phi(n)b(n,g)q^{n}$. Then
for any $\gamma = \left(%
\begin{array}{cc}
  a & b \\
  c & d \\
\end{array}%
\right) \in \Gamma_{0}(N^{2}p^{2m})$, we have the following
transformation formula,
\[
g(\phi)|_{k}\gamma = \omega(d)g(\phi_{a})
\]
where $\phi_{a}(y) = \phi(a^{-2}y)$ for all $y \in
Y_{m}=$\textbf{Z}$/Np^{m}$\textbf{Z}.
\end{lemma}
\paragraph{Proof} See \cite{Hida2} (p. 190) \qed
\paragraph{}

By a result of Hida in \cite{Hida4} (p. 24 corollary 2.3) it follows
that actually the measure $\mu_{f}$, on $\Z_{p}^{\times}$, is
cuspidal.

\paragraph{Eisenstein measure and convolution:} Of particular importance for us is the existence, which follows from
\cite{Katz1}, of the following arithmetic measure of weight one,
$dE:C(Z_{L};\ri_{K})\rightarrow \bar{S}(L;\ri_{K})$ defined by,
\[
2\int_{Z_{L}}\phi(z)dE =
\sum_{\stackrel{n=1}{(n,p)=1}}^{\infty}\left(\sum_{\stackrel{d|n}{(d,L)=1}}
sgn(d)\phi(d)\right)q^{n} \in \ri_{K}[[q]]
\]
We call this the Eisenstein-Katz measure. For a general arithmetic
measure $\mu_{g}$ of $\Z_{p}^{\times}$ associated to a modular form
of weight $\ell$ and character $\psi$ we can define a convolution
operation, see for example \cite{Hida2,PerrinRiou1}, of $\mu_{g}$
and $dE$. We consider the action of $Z_{L}$ on
$C(\Z_{p}^{\times};\ri_{K})$ by $(z \star \phi)(x):=
\psi(z)z_{p}^{\ell}\phi(z_{p}^{2}x)$ for $z \in Z_{L}$ and $\phi \in
C(\Z_{p}^{\times};\ri_{K})$. For a given integer $k \geq \ell$ and a
finite order character $\chi:Z_{L}\rightarrow \C^{\times}$ we define
the arithmetic measure $(\mu_{g}
* dE)_{\chi,k}:C(\Z_{p}^{\times};\ri_{K}) \rightarrow
\overline{S}(L;\ri_{K})$ as
\[
\int_{\Z_{p}^{\times}}\phi(x) (\mu_{g}^{L} * dE)_{\chi,k} :=
\int_{\Z_{p}^{\times}} \int_{Z_{L}}\chi(z)z_{p}^{k-1}(z^{-1} \star
\phi)(x) dE(z)d\mu_{g}^{L}(x)
\]

\section{$p$-adic Rankin-Selberg convolution}

Now we have collected all the needed background from the theory of
$p$-adic modular forms and measures to introduce $p$-adic
Rankin-Selberg convolution. In this section we state and prove a
simplified version, sufficient for our purposes, of a theorem of
Hida, as for example stated in \cite{Hida4} theorem 5.1.

Let $f \in S_{k}(\Gamma_{0}(N),\chi)$ be a normalized eigenform with
$|a(p,f)|_{p}=1$ and $(N,p)=1$. Write $f_{0} \in
\Gamma_{0}(Np),\chi)$ for the corresponding ordinary form. We recall
that $f_{0}=f-\frac{\chi(p)p^{k-1}}{u}f|[p]$ where $f|[p](z):=f(pz)$
and $u$ the root of $X^{2}-a(p,f)X+\chi(p)p^{k-1}=0$ which is a
$p$-adic unit. Let $g \in M_{\ell}(\Gamma_{0}(Jp^{\alpha}),\psi)$
with $(J,p)=1$ and $k > \ell$. Consider the cuspidal arithmetic
measure $\mu_{g}$ on $X:=\mathbb{Z}_{p}^{\times}$ that we can attach
to $g$ from the previous section. Now we assume that we can attach
to $f_{0}$ a linear form $\ell_{f_{0}}: \overline{S}_{k}(N;\ri_{K})
\rightarrow K$ as in the previous section. We also consider a
constant $c(f)$ such that $c(f)\ell_{f}$ takes integral values. Then
we have,

\begin{theorem}$($$p$-adic Rankin-Selberg convolution$)$ \label{thm:Hida} With notation as above, there is a measure $\mu_{f \times g}:\Z_{p}^{\times}\rightarrow \ri_{K}$
such that for any finite order character $\phi$ on
$\Z_{p}^{\times}$,
\[
\int_{\Z_{p}^{\times}} \phi d\mu_{f\times g} =
c(f)(-1)^{k}ta(p,f_{0})^{1-\beta}p^{\beta
\ell/2}p^{\frac{2-k}{2}\beta}\frac{D(f_{0},\mu_{g}(\phi)|_{\ell}\tau_{\beta},\ell)}{2^{k+\ell}\pi^{\ell+1}\imath^{k+\ell}<f_{0}^{\rho}\mid_{k}\tau_{Np},f_{0}>_{Np}}
\]
where,
\[
t=(-1)^{k}l.c.m(N,J)N^{k/2}J^{\ell/2}\Gamma(\ell)
\]
and $\beta$ is such that $\mu(\phi) \in
M_{\ell}(\Gamma_{1}(Jp^{\beta}))$ and $\tau_{\beta}=\left(%
\begin{array}{cc}
  0 & -1 \\
  Jp^{\beta} & 0 \\
\end{array}%
\right)$

\end{theorem}

\begin{rmk}
This is a special case of a more general result of Hida. First of
all using
    Shimura's differential operators he can show that the above p-adic measure
    interpolates the rest of the critical values of the
    Rankin-Selberg L-function i.e. $D(\ell +
    m,f_{0},\mu(\phi)|_{\ell}\tau_{\beta})$ for $0 \leq m \leq k-\ell$.
    Second, and most important, Hida can construct $p$-adic measures
    that interpolate families of modular forms (usually called
    $\Lambda$-adic modular forms), in both variables of the
    Rankin-Selberg product (under some ordinarity assumptions also
    on  the second variable).
\end{rmk}
We give the proof of the above theorem following Hida as in
\cite{Hida4} (page 76). The explicit construction of the measure
$\mu_{f \times g}$ is important for our purposes.
\paragraph{Proof:}
Let us denote by $L$ the least common multiple of $N$ and $J$. We
consider the Eisenstein-Katz measure $dE:C(Z_{L};\ri_{K})\rightarrow
\overline{S}(L;\ri_{K})$ that we have introduced in the previous
section. Recall that is defined as,
\[
2\int_{Z_{L}}\phi(z)dE =
\sum_{\stackrel{n=1}{(n,p)=1}}^{\infty}\left(\sum_{\stackrel{d|n}{(d,L)=1}}
sgn(d)\phi(d)\right)q^{n} \in \ri_{K}[[q]]
\]
and $Z_{L} = \Z_{p}^{\times} \times (\Z/L\Z)^{\times}$. We also
modify the arithmetic measure $\mu_{g}$ by defining a new one,
$\mu^{L}_{g}(\phi) := \mu_{g}(\phi)\mid [L/J]$ where
$[L/J]:\overline{S}(J;O_{K})\rightarrow \overline{S}(L;O_{K})$, as
$[L/J](\sum_{n \geq 1}a(n,g)q^{n}) \mapsto \sum_{n \geq 1}
a(n,g)q^{n\frac{L}{J}}$ , and hence $\mu^{L}_{g}$ is again an
arithmetic measure of weight $\ell$ and character $\psi$. Recall
that we have defined an action of $Z_{L}$ on
$\mathbb{Z}_{p}^{\times}$ as $z \cdot x \mapsto z_{p}^{2}x$ and by
the previous section we can consider the convoluted measure
$(\mu^{L}_{g}* dE)_{\chi,k}$, which we recall is defined by,
\[
\int_{\Z_{p}^{\times}}\phi(x) (\mu_{g}^{L} * dE)_{\chi,k} :=
\int_{\Z_{p}^{\times}} \int_{Z_{L}}\chi(z)z_{p}^{k-1}(z^{-1} \star
\phi)(x) dE(z)d\mu_{g}^{L}(x)
\]
Now we define the measure $\mu_{f\times g}$ as,
\[
\int_{\Z_{p}^{\times}} \phi d\mu_{f \times g} := c(f) \circ
\ell_{f_{0}} \circ Tr_{L/N} \circ e(\int_{\Z_{p}^{\times}}\phi
(\mu_{g}^{L} * dE)_{\chi,k})
\]
Here $Tr_{L/N}:\overline{M}(L;\ri_{K}) \rightarrow
\overline{N}(N;\ri_{K})$ is the trace operator, see
\cite{PerrinRiou1}. We do not need to give its detailed definition
but just mention that when restricted to the classical modular forms
satisfy the usual property; for $f \in S_{k}(\Gamma_{1}(N);\ri_{K})$
and $g \in M_{k}(\Gamma_{1}(L));\ri_{K})$ we have that,
$<f,Tr_{L/N}g>_{N}=<(L/N)^{k}f|[L/N],g>_{L}$. Let now $\phi$ be a
finite order character on $\Z_{p}^{\times}$. We compute the value of
our measure on $\phi$. We have,
\[
\int_{\Z_{p}^{\times}}\phi (\mu_{g}^{L} * dE)_{\chi,k} =
\int_{\Z_{p}^{\times}} \int_{Z_{L}}\chi(z)z_{p}^{k-1}(z^{-1} \star
\phi)(x) dE(z)d\mu_{g}^{L}(x) =
\]
\[
=\int_{\Z_{p}^{\times}} \int_{Z_{L}}\chi(z)z_{p}^{k-1}
\psi(z)^{-1}z_{p}^{-\ell}\phi(z_{p})^{-2}\phi(x)
dE(z)d\mu_{g}^{L}(x)=
\]
\[
=\int_{\Z_{p}^{\times}}
\int_{Z_{L}}\chi(z)z_{p}^{k-\ell-1}\psi(z)^{-1}\phi(z_{p})^{-2}\phi(x)
dE(z)d\mu_{g}^{L}=
\]
\[
=\left(\int_{\Z_{p}^{\times}}\phi(x)d\mu_{g}^{L}\right)\cdot
\left(\int_{Z_{L}}\chi\psi^{-1}(z)\phi^{-2}(z_{p})z_{p}^{k-\ell-1}dE(z)\right)
\]
Evaluating the Eisenstein measure we get,
\[
\int_{Z_{L}}\chi\psi^{-1}(z)\phi^{-2}(z_{p})z_{p}^{k-\ell-1}dE =
E_{k-\ell,Lp}(\chi\psi^{-1}\phi_{p})\mid \imath_{p}
\]
where $\phi_{p}(z)=\phi^{-2}(z_{p})$ and,
\[
E_{m,M}(\theta) :=
\frac{1}{2}L_{M}(1-m,\theta)+\sum_{n=1}^{\infty}\left(\sum_{\stackrel{0<d|n}{(d,Mp)=1}}\theta(d)d^{m-1}\right)q^{n}
\]
an Eisenstein series in $M_{m}(\Gamma_{0}(M),\theta)$. We consider
now the projection to the ordinary part. By the property $e(f\cdot
g|\imath_{p})=e(f|\imath_{p}\cdot g)$ (see \cite{Hida4}, p.24) and
since $\mu_{g}^{L}(\phi)|\imath_{p}=\mu_{g}^{L}(\phi)$ as it is
measure over $\Z_{p}^{\times}$ we have that,
\[
e(\int_{\Z_{p}^{\times}}\phi (\mu_{g}^{L} * dE)_{\chi,k})=
\left(\int_{\Z_{p}^{\times}}\phi(x)d\mu_{g}^{L}\right)\cdot
E_{k-\ell,Lp}(\chi\psi^{-1}\phi_{p})
\]
Applying the explicit formula for the linear form $\ell_{f_{0}}$ and
after writing $h := (\int_{\Z_{p}^{\times}}\phi (\mu_{g}^{L} *
dE)_{\chi,k} \in S_{k}(\Gamma_{0}(Np^{\beta})),\chi)$ we have,
\[
\int_{\Z_{p}^{\times}} \phi d\mu_{f \times g} =
c(f)a(p,f_{0})^{1-\beta}p^{(\beta-1)(k/2)}\frac{<(f_{0}^{\rho}\mid_{k}\tau_{Np})\mid
[p^{\beta-1}],Tr_{L/N}(h)>_{Np^{\beta}}}{<f_{0}^{\rho}\mid_{k}\tau_{Np},f_{0}>_{Np}}
\]
We claim the equality,
\[
<(f_{0}^{\rho}\mid_{k}\tau_{Np})\mid
[p^{\beta-1}],Tr_{L/N}(h)>_{Np^{\beta}}=(L/N)^{k/2}p^{k(1-\beta)/2}<f_{0}^{\rho}\mid\tau_{Lp^{\beta}},h>_{Lp^{\beta}}
\]
Indeed by the property of the trace operator that we described above
we have,
\begin{center}
$<(f_{0}^{\rho}\mid_{k}\tau_{Np})\mid
[p^{\beta-1}],Tr_{L/N}(h)>_{Np^{\beta}}
=(L/N)^{k}<(f_{0}^{\rho}\mid_{k}\tau_{Np})\mid [Lp^{\beta
-1}/N],h>_{Lp^{\beta}}$
\\
$=(L/N)^{k/2}p^{k(1-\beta)/2}<f_{0}^{\rho}\mid\tau_{Lp^{\beta}},h>_{Lp^{\beta}}$
\end{center}
Hence the evaluation of the measure now reads as,
\[
\int_{\Z_{p}^{\times}} \phi d\mu_{f \times g} =
c(f)a(p,f_{0})^{1-\beta}p^{(\beta-1)(k/2)}\frac{(L/N)^{k/2}p^{k(1-\beta)/2}<f_{0}^{\rho}\mid\tau_{Lp^{\beta}},h>_{Lp^{\beta}}}{<f_{0}^{\rho}\mid_{k}\tau_{Np},f_{0}>_{Np}}
\]
We note that we can write ,
\[
\mu_{g}^{L}(\phi)=(-1)^{\ell}(L/J)^{-\ell/2}(\mu_{g}(\phi)|_{\ell}\tau_{Jp^{\beta}})|_{\ell}\tau_{Lp^{\beta}}
\]
The following proposition, which is taken from \cite{Hida4} p. 63
allow us to conclude the proof.
\begin{proposition} Let $h_{1} \in S_{k}(\Gamma_{0}(Lp^{\beta}),\psi)$ and
$h_{2} \in M_{\ell}(\Gamma_{0}(Lp^{\beta}),\xi)$. Then,
\[
D(h_{1},h_{2},\ell) = t'
<h_{1}^{\rho}|_{k}\tau_{Lp^{\beta}},(h_{2}|_{\ell}\tau_{Lp^{\beta}})(E_{k-\ell,Lp}(\xi\psi))>_{Lp^{\beta}}
\]
where
$t':=2^{k+\ell}\pi^{\ell+1}(Lp^{\beta})^{\frac{1}{2}(k-\ell-2)}(\sqrt{-1})^{\ell-k}(\Gamma(\ell))^{-1}$
\end{proposition}

\section{Towards the congruences}

In this section we obtain a first form of the congruences claimed in
the introduction.
\subsection{The case $p=3$}
We start with some generalities. Let $K/\mathbb{Q}$ be a quadratic
imaginary extension of discriminant $D$ and non-trivial character
$\epsilon_{D}$. Let $\chi^{*}:\A^{\times}_{K}/K^{\times} \rightarrow
\mathbb{C}^{\times}$ be a finite order Hecke character corresponding
by class field theory to a Galois character
$\chi:Gal(K(\textrm{f}_{\chi})/K)\rightarrow \mathbb{C}^{\times}$
where $\textrm{f}_{\chi}$ the conductor of $\chi$ and
$K(\textrm{f}_{\chi})$ the ray class field for the ideal
$\textrm{f}_{\chi}$. We will also write $\chi$ for the ideal
character corresponding to $\chi^{*}$. Consider the series
$g_{\chi}(z)=\sum_{\textrm{a}\subset
\ri_{K}}\chi(\textrm{a})q^{N(\textrm{a})}$ if $\chi$ is not the
trivial character where $q = e^{2 \pi \imath z}$ and
$\chi(\textrm{a})=0$ if $(\textrm{a},\textrm{f}_{\chi})\neq 1$. In
case $\chi$ is the trivial character we define $g_{1}(z) =
\frac{1}{2}L(0,\epsilon_{D})+\sum_{\textrm{a} \subset
\ri_{K}}q^{N(\textrm{a})}$. By automorphic induction we have that
$g_{\chi}(z) \in M_{1}(\Gamma_{0}(|D|
N(\textrm{f}),\epsilon_{D}\chi|_{\mathbb{Z}})$ where by $\chi|_{\Z}$
we mean the character obtained by restricting $\chi$ to ideals in
$\Z$. Moreover it is known that for $\chi$ non-trivial we have that
$g_{\chi}(z) \in
S_{1}(\Gamma_{0}(|D|N(\textrm{f}),\epsilon_{D}\chi|_{\mathbb{Z}})$
is a primitive form.

Let us write $p$ for the prime number 3. We consider the field
$\Q(\mu_{p})/\Q)$ and we write $\p$ for the unique prime above $p$
in it. Let us now denote by $\chi$ any of the two non-trivial
character of the cyclic cubic extension
$\mathbb{Q}(\mu_{p},\sqrt[p]{m})/\mathbb{Q}(\mu_{p})$ for $m$ a
power free integer and $(m,p)=1$. Note that $\chi \equiv 1
\mod{\p}$. We consider the induced representation
$\rho:=Ind^{K}_{\Q}(\chi)$, a two dimensional Artin representation
$\rho:S_{3} \cong Gal(\mathbb{Q}(\mu_{p},\sqrt[p]{m})/\mathbb{Q})
\rightarrow GL_{2}(\mathbb{Z})$. We write $g_{\rho}$ for the
corresponding newform obtained from the discussion above with
$g_{\rho} \in S_{1}(\Gamma_{0}(m^{2}p^{r}),\epsilon_{p} \chi|_{\Z})$
where $r=1$ if $\chi$ does not ramify at $p$ and $r=3$ if it does.
Note that actually $\chi|_{\Z}$ is the trivial character. Finally
let us also write $g_{\sigma}$ for the Eisenstein series $g_{1}$.
Then, $g_{\sigma} \in M_{1}(\Gamma_{0}(p),\epsilon_{p})$.

We associate $p$-adic arithmetic measures to our modular forms
$g_{\sigma}$ and $g_{\rho}$. We modify $g_{\sigma}$ and consider the
modular form $g_{\sigma_{(m)}} := g_{\sigma}\mid \imath_{m} \in
M_{1}(\Gamma_{0}(m^{2}p,\epsilon_{p})$, with $\imath_{m}$ the
trivial character modulo $m$, i.e we remove the ``Euler factors'' at
the primes dividing $m$. We now consider the associated arithmetic
measures on $\Z_{p}^{\times}$. For $\phi \in
C(\Z_{p}^{\times};\Z_{p})$ we have,
\[
d\mu_{\rho} : \phi \mapsto \sum_{n=1}^{\infty}
\phi(n)a(n,g_{\rho})q^{n}
\]
\[
d\mu_{\sigma_{(m)}} : \phi \mapsto \sum_{n=1}^{\infty}
\phi(n)a(n,g_{\sigma_{(m)}})q^{n}
\]
Note that by construction we have that for any $\phi \in
C(\Z_{p}^{\times},\Z_{p})$,
\[
\int_{\Z_{p}^{\times}}\phi
d\mu_{\rho} \equiv \int_{\Z_{p}^{\times}}\phi d\mu_{\sigma_{(m)}}
\mod{p}
\]
where the meaning of the congruences here is term by term i.e.
$\phi(n)a(n,g_{\rho}) \equiv \phi(n)a(n,g_{\sigma_{m}}) \mod{p}$ for
all $n$.

Let now $E/\Q$ be an elliptic curve over $\Q$ with conductor $N$ and
with good ordinary reduction at $p$. Recall that we are assuming
that $E[p]$ is an irreducible $G_{\Q}$-module and moreover the Artin
conductor of the representation $\overline{\rho}_{E,p}:G_{\Q}
\rightarrow Aut(E[p])$ is equal to $N=N_{E}$. Let us write  $f \in
S_{2}(\Gamma_{0}(N);\Q)$ for the primitive form associated to $E$.
The assumption of the good ordinary reduction at $p$ implies that
$|a(p,f)|_{p}=1$ where $f(z) = \sum_{n \geq 1}a(n,f)q^{n}$ and of
course that $(N,p)=1$. Let us now write $f_{0} \in
S_{2}(\Gamma_{0}(Np);\Q)$ for the ordinary form that we can
associate to $f$ by lemma \ref{lemma:ordinary} and $\tilde{f}_{0}
\in S_{2}(\Gamma_{0}(Npm^{2});\Q)$ for the normalized eigenform that
we obtain after removing the Euler factors at $q|m$, that is
$\tilde{f}_{0}=f_{0}|\imath_{m}$. We now consider the map
$h_{2}(\Gamma_{0}(Nm^{2}p);\Z_{p}) \rightarrow \Z_{p}$ induced by
$T(n) \mapsto a(n,\tilde{f}_{0})$. Later we will prove that actually
this map, under our assumptions, induces a decomposition,
\[
h_{2}(\Gamma_{0}(Nm^{2}p);\Q_{p}) = \Q_{p} \oplus A
\]
Let us write $1_{\tilde{f}_{0}}$ for the idempotent attached to the
first summand. Moreover we consider a constant $c(f,m) \in \Z_{p}$,
defined up to $p$-adic units, such that $c(f,m)1_{\tilde{f}_{0}} \in
h_{2}(\Gamma_{0}(Nm^{2}p;\Z_{p})$. Let us now write $L$ for
$Nm^{2}$. We denote by $dE_{2,id}=dE:C(Z_{L};\Z_{p})\rightarrow
\overline{S}(L;\Z_{p})$ the Eisenstein-Katz measure on $Z_{L}$.
Recall also that for any arithmetic measure $d\mu :
C(\Z_{p}^{\times};\Z_{p}) \rightarrow \overline{S}(M;\Z_{p})$ with
$M \mid L$ we have defined another arithmetic measure $d\mu^{L}$
with values in $\overline{S}(L,\Z_{p})$ by applying the operator
$[L/M]$.

\begin{lemma} \label{lemma:congruences} Let $\phi \in C(\Z_{p}^{\times};\ri_{K}^{\times})$ be a
character of finite order. Consider the measures
$d\mu^{L}_{g_{\chi}}\ast dE$ and $d\mu^{L}_{\sigma_{(m)}} \ast dE$.
Then we have,
\[
\left|\int_{\Z_{p}^{\times}}\phi\,\,(d\mu^{L}_{g_{\rho}}\ast dE) -
\int_{\Z_{p}^{\times}}\phi\,\,(d\mu^{L}_{\sigma_{(m)}}\ast
dE)\right|_{p} < 1
\]
\end{lemma}

\paragraph{Proof} By the calculations we did in the previous section for the proof of Hida's $p$-adic Rankin-Selberg
theorem we have,
\[
\int_{\Z_{p}^{\times}}\phi\,\,(d\mu^{L}_{g_{\sigma}}\ast dE) =
(\int_{\Z_{p}^{\times}}\phi
d\mu^{L}_{g_{\rho}})(\int_{Z_{L}}\epsilon_{p}(z)\phi_{p}(z)dE)
\]
and similarly,
\[
\int_{\Z_{p}^{\times}}\phi\,\,(d\mu^{L}_{g_{\sigma_{(m)}}}\ast dE) =
(\int_{\Z_{p}^{\times}}\phi
d\mu^{L}_{g_{\sigma_{(m)}}})(\int_{Z_{L}}\epsilon_{p}(z)\phi_{p}(z)dE)
\]
with $\phi_{p}(z) = \phi^{-2}(z_{p})$. The lemma now follows from
the facts that $dE$ is an integral measure and
$|\mu^{L}_{g_{\sigma}}(\phi)-
\mu^{L}_{g_{\sigma_{(m)}}}(\phi)|_{p}<1$ as the operator $[N]$
preserves congruences. \qed

Now we are ready to prove a first type of congruences. Let us write
$u$ for $a(p,f_{0})$ and define $w$ by $uw = p$. We have,

\begin{theorem}\label{theorem:Congruences I} Consider the
quantities,
\[
R(\rho) :=
c(f,m)\alpha(\rho)\frac{P_{p}(\rho,u^{-1})}{P_{p}(\rho,w^{-1})}\frac{D_{\{p,q|m\}}(f,g_{\rho},1)}{\pi^{2}i<\tilde{f}_{0}|\tau_{Lp},\tilde{f}_{0}>_{Lp}}
\]
and
\[
R(\sigma):=c(f,m)\alpha(\sigma)\frac{P_{p}(\sigma,u^{-1})}{P_{p}(\sigma,w^{-1})}\frac{D_{\{p,q|m\}}(f,g_{\sigma},1)}{\pi^{2}i<\tilde{f}_{0}|\tau_{Lp},\tilde{f}_{0}>_{Lp}}
\]
where, $\alpha(\rho) := e_{p}(\rho)u^{-v_{p}(N_{\rho})}$,
$\alpha(\sigma):=e_{p}(\sigma)u^{-v_{p}(N_{\sigma})}$ with
 $\sigma := 1 \oplus \epsilon_{p}$
the Artin representation induced by the trivial character,
$e_{p}(\cdot)$ local epsilon factor and $v_{p}(N_{\rho})$ the p-adic
valuation of the conductor of the Artin representation. Then with
the assumptions as above $R(\rho)$ and $R(\sigma)$ are $p$-adic
integers and,
\[
R(\rho)\equiv R(\sigma) \mod{p}.
\]
\end{theorem}
Here we would like to remind the reader that under the assumption of
the elliptic curve being semi-stable we have that $D(s,f,g_{\rho}) =
L(E_{f},\rho,s)$.

\paragraph{Proof} We claim that for any character $\phi \in
C(\Z_{p}^{\times};\ri^{\times})$ we have that,
\[
|\int_{\Z_{p}^{\times}}\phi d\mu_{\tilde{f}_{0},g_{\rho}} -
\int_{\Z_{p}^{\times}}\phi
d\mu_{\tilde{f}_{0},g_{\sigma_{(m)}}}|_{p} < 1
\]
By definition we have,
\[
\int_{\Z_{p}^{\times}}\phi d\mu_{\tilde{f}_{0},g_{\rho}} = c(f,m)
\ell_{\tilde{f}_{0}} \circ
e(\int_{\Z_{p}^{\times}}\phi\,\,(d\mu^{L}_{g_{\rho}}\ast dE))
\]
\[
\int_{\Z_{p}^{\times}}\phi d\mu_{\tilde{f}_{0},g_{\sigma}} = c(f,m)
\ell_{\tilde{f}_{0}} \circ
e(\int_{\Z_{p}^{\times}}\phi\,\,(d\mu^{L}_{g_{\sigma}}\ast dE))
\]
Note that the trace operator is now just the identity. By the
definition of the linear form $\ell_{\tilde{f}_{0}}$ we have that
$c(f,m) \ell_{\tilde{f}_{0}}(e\, h) =
a(1,c(f,m)1_{\tilde{f}_{0}}e\,h)$ for $h \in
\overline{S}(L;\Z_{p})$. But we have that $|e\,h|_{p} \leq |h|_{p}$
and also $|c(f,m)1_{\tilde{f}_{0}}\,e\,h|_{p} \leq |e\,h|_{p}$ and
hence by lemma \ref{lemma:congruences} we establish the claim. Now
in order to obtain the congruences we evaluate both measures at the
trivial character $\imath_{p}$ modulo $p$ and hence we have,
\[
\int_{\Z_{p}^{\times}}\imath_{p} d\mu_{\tilde{f}_{0},g_{\rho}}
\equiv \int_{\Z_{p}^{\times}}\imath_{p}
d\mu_{\tilde{f}_{0},g_{\sigma_{(m)}}} \mod{p}
\]
We now work both sides of the above equation. We start with the left
hand side. By theorem ~\ref{thm:Hida} we have,
\[
\int_{\Z_{p}^{\times}}\imath_{p} d\mu_{\tilde{f}_{0},g_{\rho}}=
c(f,m) utp^{\beta/2} u^{-\beta}
\frac{D(\tilde{f}_{0},\mu_{g_{\rho}}(\imath_{p},1)|_{1}\tau_{\beta})}{\pi^{2}i<\tilde{f}_{0}|_{2}\tau_{Lp},\tilde{f}_{0}>_{Lp}}
\]
where $\beta$ is such that $\mu_{g_{\rho}}(\imath_{p}) =
g_{\chi}|\imath_{p}\in S_{1}(\Gamma_{1}(m^{2}p^{\beta})$. We
consider the Rankin-Selberg product
$D(\tilde{f}_{0},\mu_{g_{\rho}}(\imath_{p})|_{1}\tau_{\beta},1) =
D(\tilde{f}_{0},g_{\rho}|\imath_{p})|_{1}\tau_{\beta},1)$. We will
write $g$ for $g_{\rho}$ and $M$ for $m^{2}$. Let us assume first
that the character $\chi$ is not ramified above $p$ and hence $\beta
= 2$. We can write in this case $g|\imath_{p} = g -
a(p,g)g|[p]=g-g|[p]$ as $a(p,g)=1$. We apply $\tau_{Mp^{2}} = \left(%
\begin{array}{cc}
  0 & -1 \\
  Mp^{2} & 0 \\
\end{array}%
\right)$ to the above equation and we use the fact that $g \in
\Gamma_{1}(Mp)$ is a primitive form of level $Mp$ and hence satisfies $g|_{1}\left(%
\begin{array}{cc}
  0 & -1 \\
  Mp & 0 \\
\end{array}%
\right) = W(g)g^{\rho}=W(g)g$, as $g$ has rational coefficients. The
quantity $W(g)$ is usually called the root number of $g$. We have,
\[
(g|\imath_{p})|_{1}\tau_{Mp^{2}}=g|_{1}\tau_{Mp^{2}}-g|[p]|_{1}\tau_{Mp^{2}}=
\]
\[
=g|_{1}\left(%
\begin{array}{cc}
  0 & -1 \\
  Mp & 0 \\
\end{array}%
\right)\left(%
\begin{array}{cc}
  p & 0 \\
  0 & 1 \\
\end{array}%
\right)-p^{-\frac{1}{2}}g|_{1}\left(%
\begin{array}{cc}
  p & 0 \\
  0 & 1 \\
\end{array}%
\right)\left(%
\begin{array}{cc}
  0 & -1 \\
  Mp^{2} & 0 \\
\end{array}%
\right)=
\]
\[
= W(g)g|_{1}\left(%
\begin{array}{cc}
  p & 0 \\
  0 & 1 \\
\end{array}%
\right)-p^{-\frac{1}{2}}g|_{1}\left(%
\begin{array}{cc}
  0 & -1 \\
  Mp & 0 \\
\end{array}%
\right)\left(%
\begin{array}{cc}
  p & 0 \\
  0 & p \\
\end{array}%
\right)=
\]
\[
=p^{\frac{1}{2}}W(g)g|[p]-W(g)p^{-\frac{1}{2}}g =
-p^{-\frac{1}{2}}W(g)(g-pg|[p])
\]
Hence we get that,
\[
D(f_{0},(g|\imath_{p})|_{1}\tau_{Mp^{2}},1)=
-p^{-\frac{1}{2}}W(g)(1-a(p,f_{0}))D(f_{0},g,1)=p^{-\frac{1}{2}}W(g)u
P_{p}(g,u^{-1})D(f_{0},g,1)
\]
Moreover we note that $D(\tilde{f}_{0},g,1) =
P_{p}(g,w^{-1})^{-1}D_{\{p,q|m\}}(f,g,1)$ where we have removed the
Euler factor at $p$ and $q|m$ from the primitive Rankin-Selberg
product $D(f,g,1)$. Also, Balister in \cite{Balister} p.17 has
computed the local epsilon factors of $\rho$ from where we get
$e_{p}(\rho) = p^{\frac{1}{2}}W_{p}(g)$ and $W_{q}(g) = 1$ for $q|m$
and hence we conclude,
\[
\int_{\Z_{p}^{\times}}\imath_{p}d\mu_{\tilde{f}_{0},g_{\rho}}=
c(f,m) ute_{p}(\rho) u^{-1}
\frac{P_{p}(\rho,u^{-1})}{P_{p}(\rho,w^{-1})}\frac{D_{\{p,q|m\}}(f,g,1)}{\pi^{2}i<\tilde{f}_{0}|_{2}\tau_{Lp},\tilde{f}_{0}>_{Lp}}
\in \Z_{p}
\]
The case where $\chi$ is ramified at $p$ is easier as
$g|\imath_{p}=g$ since $P_{p}(\rho,X)=1$ and hence we can use
directly the action of $\tau_{Mp^{3}}$. Hence also in this case we
get,
\[
\int_{\Z_{p}^{\times}}\imath_{p}d\mu_{\tilde{f}_{0},g_{\rho}}= c(f)
ute_{p}(\rho) u^{-3}
\frac{D_{\{p,q|m\}}(f,g,1)}{\pi^{2}i<\tilde{f}_{0}|_{2}\tau_{Lp},\tilde{f}_{0}>_{Lp}}
\in \Z_{p}
\]
We now work the right hand side of the congruences. We have that,
\[
E_{1}(\epsilon_{p})(z) := g_{\sigma}(z) =
\frac{1}{2}L(0,\epsilon_{p}) +
\sum_{n=1}^{\infty}(\sum_{0<d|n}\epsilon_{p}(d))q^{n} \in
M_{1}(\Gamma_{0}(p),\epsilon_{p})
\]
where we write, as always, $\epsilon_{p}$ for the non-trivial
character of $Gal(\Q(\mu_{3})/\Q)$. We consider now the imprimitive
Rankin-Selberg $L$-function
$D(\tilde{f}_{0},g_{1}|\imath_{m}|\imath_{p}|_{1}\tau_{m^{2}p^{2}},1)$
where $\tau_{m^{2}p^{2}} = \left(%
\begin{array}{cc}
  0 & -1 \\
  m^{2}p^{2} & 0 \\
\end{array}%
\right)$. We can consider each prime separately, i.e. first we
consider the quantity $D(\tilde{f}_{0},g_{1}|\imath_{p}|\tau_{p},1)$
with
$\tau_{p} = \left(%
\begin{array}{cc}
  0 & -1 \\
  p^{2} & 0 \\
\end{array}%
\right)$ and then the quantity
$D(\tilde{f}_{0},g_{1}|\imath_{q}|\tau_{q},1)$
with $\tau_{q} = \left(%
\begin{array}{cc}
  0 & -1 \\
  pq^{2} & 0 \\
\end{array}%
\right)$ for $q|m$. Let assume this and do the calculations and at
the end we return to this point. As before we can write
$E_{1}(\epsilon_{p})\mid \imath_{p} =
E_{1}(\epsilon_{p})-E_{1}(\epsilon_{p})|_{1}[p]$. Working as above
we have,
\[
E_{1}(\epsilon_{p})\mid \imath_{p} \mid \left(%
\begin{array}{cc}
  0 & -1 \\
  p^{2} & 0 \\
\end{array}%
\right) =
p^{-\frac{1}{2}}W(E_{1}(\epsilon_{p}))(E_{1}(\epsilon_{p})-pE_{1}(\epsilon_{p})\mid
[p])
\]
Hence we have,
\[
D(\tilde{f}_{0},E_{1}(\epsilon_{p})| \imath_{p} |_{1} \tau_{p},1)=
p^{-\frac{1}{2}}W(E(\epsilon_{p}))a(p,\tilde{f}_{0})(1-a(p,\tilde{f}_{0})^{-1})D(\tilde{f}_{0},E_{1}(\epsilon_{p}),1)
\]
But for the Eisenstein series $E_{1}(\epsilon_{p})$ we know that
$D(\tilde{f}_{0},E_{1}(\epsilon_{p}),1)=L(\tilde{f}_{0},1)L(\tilde{f}_{0},\epsilon_{p},1)$
or equivalently $D(\tilde{f}_{0},E_{1}(\epsilon_{p}),1) =
(1-a(p,\tilde{f}_{0})p^{-1})^{-1}L_{\{p\}}(\tilde{f},1)L_{\{p\}}(\tilde{f},\epsilon_{p},1)$.
Recall that we have defined $u:=a(p,\tilde{f}_{0})=a(p,f_{0})$ and
$uw=p$ we finally get,
\[
D(\tilde{f}_{0},E_{1}(\epsilon_{p},1)| \imath_{p}|_{1} \tau_{p})=
p^{-\frac{1}{2}}W(E(\epsilon_{p}))u\frac{1-u^{-1}}{1-w^{-1}}L_{\{p,q|m\}}(f,1)L_{\{p,q|m\}}(f,\epsilon_{p},1)
\]
We now consider a prime $q\mid m$. Now we write just $g$ for
$g_{1}=E_{1}(\epsilon_{p})$. As $q \neq p$, we have $g \mid
\imath_{q} = g - a(q,g)g\mid[q]+\epsilon_{p}(q)g \mid [q^{2}]$. Now we apply the operator $\tau_{q} = \left(%
\begin{array}{cc}
  0 & -1 \\
  q^{2}p & 0 \\
\end{array}%
\right)$. Doing the calculations as before we get
$g|\imath_{q}\mid_{1}\tau_{q} =
q^{-1}\epsilon_{p}(q)W(g)(g-\epsilon_{p}(q)a(q,g)qg\mid[q]+\epsilon_{p}(q)q^{2}g\mid[q^{2}])$.
But note that since $L(s,\tilde{f}_{0})$ has no Euler factors at $q
\mid m$ we have that $D(s,\tilde{f}_{0},g|[q^{r}])=0$ for any $r
\geq 1$. So we obtain $D(\tilde{f}_{0},g|\imath_{q}\mid \tau_{q},1)
= q^{-1}\epsilon_{p}(q)W(g)D(\tilde{f}_{0},g,1)$. Putting all
together and noticing that $q^{-1}\epsilon_{p}(q) \equiv 1 \mod{p}$
and $W_{q}(g)=1$ we get,
\[
\int_{\Z_{p}^{\times}}\imath_{p}d\mu_{f,g_{\sigma_{m}}}= c(f,m)
ute_{p}(\sigma) u^{-1}
\frac{P_{p}(\sigma,u^{-1})}{P_{p}(\sigma,w^{-1})}\frac{D_{\{p,q|m\}}(f,g_{\sigma},1)}{\pi^{2}i<\tilde{f}_{0}|_{2}\tau_{Lp},\tilde{f}_{0}>_{Lp}}
\in \Z_{p}
\]
The fact that $|ut|_{p}=1$ allows us to conclude the proof of the
theorem. Let us now also justify our claim that we can work each
prime separately. For simplicity we do the case of $m=q$ but we will
become obvious how one obtains the general case. So with $g$ as
above we have,
\[
g|i_{pq} = (g - g|[p])-a(q,g)((g - g|[p]))|[q]+\epsilon_{p}(q)((g -
g|[p]))[q^{2}]
\]
Now we apply the operator $\tau_{p^{2}q^{2}}$. We claim that only
the term $\epsilon_{p}(q)((g - g|[p]))[q^{2}]|_{1}\tau_{p^{2}q^{2}}$
will survive after considering the Rankin-Selberg convolution with
$\tilde{f}_{0}$. Indeed as $\tilde{f}_{0}$ has no Euler factors at
$m$ its Rankin-Selberg convolution with a form $g'$ with $a(n,g')=0$
if $(n,q)=1$, will be trivial. Consider now,
$g|[q^{i}p^{j}]|\tau_{p^{2}q^{2}}$ with $i=0,1,2$ and $j=0,1$. Then,
\[
g|[q^{i}p^{j}]|\tau_{p^{2}q^{2}}=\frac{1}{(q^{i}p^{j})^{\frac{1}{2}}}g|_{1}\left(%
\begin{array}{cc}
  q^{i}p^{j} & 0 \\
  0 & 1 \\
\end{array}%
\right)\left(%
\begin{array}{cc}
  0 & -1 \\
  p^{2}q^{2} & 0 \\
\end{array}%
\right)=
\]
\[
=\frac{1}{(q^{i}p^{j})^{\frac{1}{2}}}g|_{1}\left(%
\begin{array}{cc}
  0 & -1 \\
  p & 0 \\
\end{array}%
\right)\left(%
\begin{array}{cc}
  pq^{2} & 0 \\
  0 & q^{i}p^{j} \\
\end{array}%
\right)=\frac{1}{(q^{i}p^{j})^{\frac{1}{2}}}W(g)g|_{1}\left(%
\begin{array}{cc}
  q^{i}p^{j} & 0 \\
  0 & q^{i}p^{j} \\
\end{array}%
\right)\left(%
\begin{array}{cc}
  p^{1-j}q^{2-i} & 0 \\
  0 & 1 \\
\end{array}%
\right)=
\]
\[
=\frac{1}{(q^{i}p^{j})^{\frac{1}{2}}}W(g)g|_{1}\left(%
\begin{array}{cc}
  p^{1-j}q^{2-i} & 0 \\
  0 & 1 \\
\end{array}%
\right)=(\frac{p^{1-j}q^{2-i}}{(q^{i}p^{j})})^{\frac{1}{2}}W(g)g|[p^{1-j}q^{2-i}]
\]
So we see that if $i\neq 2$ then
$D(1,\tilde{f}_{0},g|[q^{i}p^{j}]|\tau_{p^{2}q^{2}})=0$. Moreover we
see from the above computations that the term $\epsilon_{p}(q)((g -
g|[p]))[q^{2}]$ equals
$\frac{\epsilon_{p}(q)}{q}p^{-\frac{1}{2}}W(g)(g-pg|[p])$, which
concludes our claim. Now it is not hard to check that our argument
extends to the general case. One has again to observe that only
terms of the form $g'|[m^{2}]$ will survive after the Rankin-Selberg
convolution. \qed

\subsection{The case $p>3$}

Let us fix some notation first. Let us write $K:=\Q(\mu_{p})$ and
$F:=\Q(\mu_{p})^{+}$, a totally real field and $[K:F]=2$. Moreover
we write $\p$ for the unique prime above $p$ in $F$ and
$\epsilon_{\p}$ for the non-trivial character of $K/F$. Let us also
denote by $\chi$ a non-trivial character of the cyclic extension
$K(\sqrt[p]{m})/K$ of degree $p$ for some $p^{th}$ power free
integer $m$. Moreover we write $\rho := Ind_{K}^{F}(\chi)$ and $R :=
Ind_{F}^{\Q}(\rho)=Ind_{K}^{\Q}(\chi)$. Also we write $\sigma :=
Ind_{K}^{F}(1)=1 \oplus \epsilon_{\p}$ and $\Sigma=
Ind_{K}^{\Q}(1)=\oplus_{r=1}^{p}\theta^{r}$ for some character
$\theta$ of $Gal(K/\Q)$. Our aim is to establish congruences between
the quantities
\[
Q(R) := e_{p}(R)u^{-v_{p}(N_{R})}
\frac{P_{p}(R,u^{-1})}{P_{p}(R,w^{-1})}\frac{L_{S}(E/\Q,R,1)}{(\Omega(E)_{+}\Omega(E)_{-})^{\frac{p-1}{2}}}
\]
\[
Q(\Sigma) := e_{p}(\Sigma)u^{-v_{p}(N_{\Sigma})}
\frac{P_{p}(\Sigma,u^{-1})}{P_{p}(\Sigma,w^{-1})}\frac{L_{S}(E/\Q,\Sigma,1)}{(\Omega(E)_{+}\Omega(E)_{-})^{\frac{p-1}{2}}}
\]
where $S$ is the set of primes consisting of $p$ and $q|m$. From the
inductive properties of the $L$ functions we note the equalities
$L_{S}(E/\Q,R,1)=L_{S}(E/\Q,Ind_{K}^{\Q}(\chi),1)=L_{S}(E/\Q,Ind_{K}^{F}Ind_{F}^{\Q}(\chi),1)=
L_{S}(E/F,Ind_{K}^{F}(\chi),1)=L_{S}((E/F,\rho,1)$ and also in the
same way $L_{S}(E/\Q,\Sigma,1)=L_{S}(E/F,\sigma,1)$. Of course here
the set $S$ contains the primes of $\ri_{F}$ and $\ri_{K}$ above $m$
and of course $\p$. Moreover as the inductive properties hold for
Euler factors we can also conclude that $P_{p}(R,X)=P_{\p}(\rho,X)$
and $P_{p}(\Sigma,X) = P_{\p}(\sigma,X)$. For the local epsilon
factors and the conductor we know that they are inductive in degree
zero and so we have that,
$\frac{e_{p}(R)}{e_{p}(\Sigma)}=\frac{e_{\p}(\rho)}{e_{\p}(\sigma)}$.
We now consider the quantities,
\[
Q(\rho) := e_{\p}(\rho)u^{-v_{\p}(N_{\rho})}
\frac{P_{\p}(\rho,u^{-1})}{P_{\p}(\rho,w^{-1})}\frac{L_{S}(E/F,\rho,1)}{(\Omega(E)_{+}\Omega(E)_{-})^{\frac{p-1}{2}}}
\]
\[
Q(\sigma) := e_{\p}(\sigma)u^{-v_{\p}(N_{\sigma})}
\frac{P_{\p}(\sigma,u^{-1})}{P_{\p}(\sigma,w^{-1})}\frac{L_{S}(E/F,\sigma,1)}{(\Omega(E)_{+}\Omega(E)_{-})^{\frac{p-1}{2}}}
\]
Let us now write $f \in S_{2}(\Gamma_{0}(N),\Z)$ for the primitive
cusp form that we can associate to $E$ and by $\pi_{f}$ the
corresponding cuspidal automorphic representation, i.e.
$L(E/\Q,s)=L(f,s) = L(\pi_{f},s)$. Notice that as $F/\Q$ is cyclic
we can consider the base change of $\pi_{f}$ from $\Q$ to $F$, a
cuspidal automorphic representation $\pi_{\phi}$ of $GL(2,\A_{F})$
such that $L(\pi_{\phi},s)=\prod_{r=1}^{\frac{p-1}{2}}L(\pi_{f}
\otimes \eta^{r},s)$ for a finite order Hecke character $\eta$ that
corresponds to a Galois character that generates $Gal(F/\Q)^{\vee}$.
Let us write $\phi$ for the Hilbert modular form of parallel weight
two that we attach to $\pi_{\phi}$ in the canonical way that we have
described in the introduction.

By automorphic induction for degree two extensions we can associate
a Hilbert modular form $\g_{\rho} \in
S_{1}(N_{K/F}(\f_{\chi})\p,\epsilon_{\p})$ to $\rho$ (the character
is just $\epsilon_{\p}$ as $\chi$ is anti-cyclotomic), and an
Eisenstein series $\E_{\sigma} \in M_{1}(\p,\epsilon_{\p})$ to
$\sigma$. From now on we will write $M$ for the ideal
$N_{K/F}(\f_{\chi})$ of $\ri_{F}$. Note that as we assume that $E$
is semi-stable we have that $L(E/F,\chi,1)=D(\phi,g_{\rho},1)$ and
$L(E/F,\sigma,1)=D(\phi,\E_{\sigma},1)$. Moreover we write
$\phi_{0}$ for the ``ordinary'' Hilbert modular form that one can
attach to $\phi$ such that $C(\q,\phi_{0})=C(\q,\phi)$ if $\q\neq
\p$ and $C(\p,\phi_{0})$ is the $p$-adic unit root of the equation
$x^{2}-C(\p,\phi)+p=0$. Note that $C(\p,\phi)=a(p,f)$. Finally by
$\tilde{\phi}_{0}$ we denote the Hilbert modular form that we obtain
from $\phi_{0}$ by removing the Euler factors above $m$.

We can extend all that we did above for the case of $p=3$ to the
more general setting of $p>3$, where now instead of working with
elliptic modular forms we work with Hilbert modular forms. In
particular we could introduce the notion of a $p$-adic Hilbert
modular form as in Hida \cite{Hida6} using the $q$-expansion
principle or their moduli interpretation as in Katz \cite{Katz2}.
However as we said in the introduction we do not have yet a theorem
for the general case for reasons that will explain later. So in this
section we restrict ourselves to just state the following theorem, a
proof of which can be found in \cite{BougThesis}. It will be enough
in order to address the issues that prevent us from proving a
general theorem for $p>3$.
\begin{theorem} Let $\gamma$ be the power of $\p$ in the level of $\E_{\sigma}$ and
$\beta$ the power of $\p$ in the level of $\g_{\rho}$. Consider the
quantities,
\[
Q(\rho):=C(\p,\phi_{0})^{-(\beta-1)}p^{\frac{\beta}{2}}
\frac{D(\tilde{\phi}_{0},\g_{\rho}|_{\imath_{\p}}|\tau_{M\p^{\beta}},1)}{i^{\frac{p-1}{2}}\pi^{p-1}<\tilde{\phi}_{0}|\tau_{L\p},\tilde{\phi}_{0}>_{L\p}}
\]
\[
Q(\sigma):=C(\p,\phi_{0})^{-(\gamma-1)}p^{\frac{\gamma}{2}}
\frac{D(\tilde{\phi}_{0},\E_{\sigma}|_{\imath_{M\p}}|\tau_{M\p^{\gamma}},1)}{i^{\frac{p-1}{2}}\pi^{p-1}<\tilde{\phi}_{0}|\tau_{L\p},\tilde{\phi}_{0}>_{L\p}}
\]
Then there exists a well-determined constant $c(\phi,m) \in
\ri_{\p}$ depending only on $\phi$ and $m$ such that, both
$c(\phi,m)Q(\rho)$ and $c(\phi,m) Q(\sigma)$ are $p$-adically
integral and,
\[
c(\phi,m)Q(\rho) \equiv c(\phi,m) Q(\sigma) \mod{\p}
\]
\end{theorem}

\section{Congruences between special values}

In this section our aim is to obtain a better understanding of the
nature of the constant $c(f,m)$ appearing in theorem
\ref{theorem:Congruences I} and its relation with the choice of our
periods. Vaguely speaking, we show that the reason that this
constant appears is that the Petersson inner product is not the
right choice to obtain integral values of the ratio
(L-values)/(automorphic periods) and this constant measures this
failure.

We start by showing that our map $h_{2}(\Gamma_{0}(Npm^{2});\Z_{p})
\rightarrow \Z_{p}$ given by $T(n) \mapsto a(n,\tilde{f}_{0})$,
under the assumptions stated in the introduction, induces a
decomposition of the form, $h_{2}(\Gamma_{0}(Npm^{2});\Q_{p}) =
\Q_{p} \times A$. This will be done by showing that actually factors
through a local ring of $h_{2}(\Gamma_{0}(Npm^{2});\Z_{p})$ that is
reduced. Then our next goal is to relate the quantity $c(f,m)$ to
the periods $<\tilde{f}_{0}|\tau_{N_{\Sigma}},\tilde{f}_{0}>$, with
$N_{\Sigma}=N_{E}m^{2}p$ that appear in the congruences of theorem
\ref{theorem:Congruences I}. Namely we will show that there is a
period determinant, we call it $\Omega(f)_{\Sigma}$ for $\Sigma$ the
set of primes dividing $m$, such that $c(f,m)\Omega(f)_{\Sigma} =
<\tilde{f}_{0}|\tau_{N_{\Sigma}},\tilde{f}_{0}>$, up to a $p$-adic
unit.

In order to conclude the theorem we will show that actually the
quantity $\Omega(f)_{\Sigma}$ is independent of $m$ and the prime
$p$. Hence we are reduced down to the primitive form $f$ and the
existence of a strong parametrization of $E$ by $X_{0}(N)$ allow us
to obtain a relation with the N\'{e}ron periods.

As we mentioned in the introduction, in this section we rely on
Wiles' deep results in \cite{Wiles}. Indeed the local ring through
which our map factors is nothing else than the reduced algebra
$\T_{\Sigma}$, using Wiles' notation, that he eventually proves to
be isomorphic to some universal deformation rings that classify
Galois representations with some predetermined properties that
deform a fixed modular mod $p$ representation. Actually for our
purposes we need the minimal level i.e. $\Sigma=\emptyset$, to be
the conductor of the elliptic curve and this is the reason for
imposing the assumption that the level of the elliptic curve is the
same with the conductor of the reduced (mod 3) representation.

Then we define the module of congruences using this reduced local
ring. It measures congruences between our modified form
$\tilde{f}_{0}$ and other normalized eigenforms. The annihilator of
this module will be eventually the quantity $c(f,m)$. We will
compare this module of congruences with what may be called the
cohomological module of congruences. Under our assumptions, Wiles'
results on the freeness of $H^{1}(X_{0}(N_{\Sigma}),\Z_{p})_{\m}$ as
a $h_{2}(\Gamma_{0}(Npm^{2});\Z_{p})_{\m}$-module, where here we
write $\m$ for the maximal ideal that corresponds to the form that
one obtains reducing $\tilde{f}_{0}$ modulo $p$, will give that
actually the two modules are isomorphic. This will relate $c(f,m)$
to the periods $<\tilde{f}_{0}|\tau_{N_{\Sigma}},\tilde{f}_{0}>$.
Then again a deep result of Wiles, the generalization of the
so-called ``Ihara's lemma'' will essentially say that this relation
does not depend on the change of level that we have introduced to
$f$ by removing the Euler factors at $m$ and modifying the one at
$p$. Hence we can reduce our study to the initial level $N$ where
the modularity of the elliptic curve provides us the way to obtain
the relation with N\'{e}ron periods.

We would like here to mention that in this section we have benefited
the most from the article of Darmon, Diamond and Taylor \cite{DDT}
based on Wiles' fundamental paper \cite{Wiles}. Most of the
constructions and proofs here are minor modifications, mainly just
restricting their constructions to our specific case, of the ones
done in their paper.

\paragraph{Structure of Hecke algebras:} In this section we collect
some well known facts about the structure of the integral Hecke
algebra $h_{2}(\Gamma_{0}(N);\Z)$ and its completion
$h_{2}(\Gamma_{0}(N);\Z_{p})$ at some prime $p$. Our main reference
is \cite{DDT}. We fix a finite extension $K$ of $\Q_{p}$ and we
denote by $\ri_{K}$ the ring of integers. We write $\la$ for the
maximal ideal in $\ri_{K}$ and $k$ for $\ri/\la$. Let us moreover
fix algebraic closures $\overline{K}$ and $\overline{k}$ of $K$ and
$k$. We consider the Hecke algebras (a)
$h_{2}(\Gamma_{0}(N);\ri_{K})= h_{2}(\Gamma_{0}(N);\Z)\otimes_{\Z}
\ri_{K}$, (b) $h_{2}(\Gamma_{0}(N);K)=h_{2}(\Gamma_{0}(N);\ri_{K})
\otimes_{\ri_{K}} K$ and (c)
$h_{2}(\Gamma_{0}(N);k)=h_{2}(\Gamma_{0}(N);\ri_{K})\otimes_{\ri_{K}}k$.
From the going-up and going-down theorems we have that the maximal
prime ideals $\m \subset h_{2}(\Gamma_{0}(N);\ri_{K})$ are above the
prime $\lambda$ i.e. $\m \cap \ri_{K} = (\la)$ and for the minimal
primes $\p$, $\p \cap \ri_{K} = (0)$. Moreover we have the
isomorphism (\cite{DDT}, p.90) $h_{2}(\Gamma_{0}(N);\ri_{K}) \isom
\prod_{\m}h_{2}(\Gamma_{0}(N);\ri_{K})_{\m}$ where the product is
over the finitely many maximal ideals of
$h_{2}(\Gamma_{0}(N);\ri_{K})$.

Let $f \in S_{2}(\Gamma_{0}(N);\overline{K})$ be a normalized
eigenform and consider the K-algebra homomorphism $\lambda_{f}:
h_{2}(\Gamma_{0}(N);K) \rightarrow \overline{K}$, that sends $T(n)
\mapsto a(n,f)$. Using this we can associate with $f$ a maximal
ideal of $h_{2}(\Gamma_{0}(N);K)$ by $\ker(\lambda_{f})$ which
depends only on the $G_{K}$ conjugacy class of $f$. In the same way
we can associate a maximal ideal of $h_{2}(\Gamma_{0}(N);k)$ to a
normalized eigenform $g \in S_{2}(\Gamma_{0}(N);\overline{k})$. We
usually write $\overline{f}$ for the reduction modulo $\lambda$ of a
form $f$ with integral Fourier expansion. The following proposition
is taken from \cite{DDT}, p. 90.

\begin{proposition} \label{proposition:structure}Let us denote by $S_{2}(N;\overline{K})(G_{K})$ the normalized eigenforms
in $S_{2}(\Gamma_{0}(N);\overline{K})$ up to $G_{K}$ conjugacy and
by $S_{2}(N;\overline{k})(G_{k})$ the normalized eigenforms in
$S_{2}(\Gamma_{0}(N);\overline{k})$ up to $G_{k}$ conjugacy. Then
the elements in $S_{2}(N;\overline{K})(G_{K})$ $($resp in
$S_{2}(N;\overline{k})(G_{k}))$ are in bijection with the maximal
ideals of $h_{2}(\Gamma_{0}(N);K)$ $($resp maximal ideals of
$h_{2}(\Gamma_{0}(N);k))$ which in turn are in bijection with the
minimal primes of $h_{2}(\Gamma_{0}(N);\ri_{K})$ $($maximal primes
of $h_{2}(\Gamma_{0}(N);\ri_{K}) \})$.
\end{proposition}

Finally we note that if we let $\m$ be a maximal ideal in
$h_{2}(\Gamma_{0}(N);\ri_{K})$ and consider the maximal ideals $\p
\subset h_{2}(\Gamma_{0}(N);K)$ with $\p \cap
h_{2}(\Gamma_{0}(N);\ri_{K}) \subset \m$ then we have an isomorphism
$h_{2}(\Gamma_{0}(N);\ri_{K})_{\m} \otimes_{\ri_{K}} K \isom
\prod_{\p} h_{2}(\Gamma_{0}(N);K)_{\p}$. We also mention what the
Atkin-Lehner theory implies for the Hecke ring
$h_{2}(\Gamma_{0}(N);K)$ under the assumption that the field $K$
contains all the coefficients of all the primitive forms of
conductor dividing $N$ and trivial character. Let us denote by
$\mathcal{P}(N)$ the set of primitive forms of conductor dividing
$N$. Then we have that, $S_{2}(\Gamma_{0}(N);K) = \oplus_{f \in
\mathcal{P}(N)}S_{K,f}$ where $S_{K,f}$ is the $K$-linear span of
$\{f(\alpha z): \alpha \mid N/N_{f}\}$ with $N_{f}$ the conductor of
the primitive form $f$. For each $f=\sum_{n \geq 0} a(n,f)q^{n}$ in
$\mathcal{P}(N)$ we denote by $h_{2}(\Gamma_{0}(N);K)[f]$ the image
of $h_{2}(\Gamma_{0}(N);K)$ in $End_{K}(S_{K,f})$. We consider the
polynomial ring, $A_{K,f}=K[u_{f,q}:\forall q \mid N/N_{f}]$ and the
ideal $I_{K,f} \subset A_{K,f}$ generated by the polynomials
$P_{f,q}(u_{f,q}) =
u_{f,q}^{v_{q}(N/N_{f})-1}(u_{f,q}^{2}-a(q,f)u_{f,q}+1_{(N_{f})}q)$
where we write $v_{q}(N/N_{f})$ for the valuation at $q$. We now
prove the following, which is a version of lemma 4.4 of \cite{DDT},
\begin{lemma}\label{lemma:Hecke algebras}
There is an isomorphism of $K$-algebras $\phi :
h_{2}(\Gamma_{0}(N);K) \isom \prod_{f \in \mathfrak{N}}
A_{K,f}/I_{K,f}$ defined by,
\[
\phi(T(q))_{f} :=
\left\{%
\begin{array}{ll}
   a(q,f), & \hbox{if $(q, N/N_{f})=1$;} \\
   u_{f,q} \mod{I_{K,f}}, & \hbox{otherwise.} \\
\end{array}%
\right.
\]
\end{lemma}

\paragraph{Proof}
We define the $K$-algebra homomorphism, $\Theta_{f} : A_{K,f}
\rightarrow h_{2}(\Gamma_{0}(N);K)[f]$ by $u_{f,q} \mapsto T(q)$.
Notice that the polynomial $P_{f,q}$ is the characteristic
polynomial of the operator $T(q)$ acting on the space spanned by the
forms $\{f(\alpha q^{i}z):i=1,\ldots, v_{q}(N/N_{f})\}$ for each
$\alpha$ dividing $N/N_{f}q^{v_{q}(N/N_{f})}$. That implies that
$I_{K,f}$ is in the kernel $\Theta_{f}$. Hence, we have the
following surjection, which we denote by $\Theta$,
\[
\Theta: \prod_{f}A_{K,f}/I_{K,f} \twoheadrightarrow
\prod_{f}h_{2}(\Gamma_{0}(N);K)[f]
\]
But since $h_{2}(\Gamma_{0}(N);K) \hookrightarrow
\prod_{f}h_{2}(\Gamma_{0}(N);K)[f]$ the following counting argument
establishes the isomorphism,
\[
dim_{K}(h_{2}(\Gamma_{0}(N);K))=dim_{K}S_{2}(\Gamma_{0}(N),K)=\sum_{f}\sigma_{0}(N/N_{f})=\sum_{f}dim_{K}A_{K,f}/I_{K,f}
\]
with $\sigma_{0}(n)=\sum_{0<d|n}1$. \qed

\paragraph{Reduced Hecke algebras:} Recall that we are interested in the $K$-algebra homomorphism
$h_{2}(\Gamma_{0}(Npm^{2}),K) \rightarrow K$ that corresponds to a
normalized eigenform that arises from an ordinary primitive form $f$
of conductor $N$ after ``removing'' the Euler factors at the primes
dividing $m$ and modifying the Euler factor at $p$ by keeping only
the $p$-adic unit part. In this section we are going to show that
under some assumptions on $f$, that we now describe, this
homomorphism factors through a local ring
$h_{2}(\Gamma_{0}(Nm^{2}p,\ri)_{\m}$ that is reduced. We recall that
if we write $\rho_{f}:=\rho_{f,p}:G_{\Q} \rightarrow GL_{2}(\Z_{p})$
for the $p$-adic representation attached to $f$, then we assume that
the reduced, mod p, representation $\overline{\rho}:=\rho \mod{p}$,
is irreducible. Let us now write $N(\bar{\rho})$ for the Artin
conductor of $\bar{\rho}$. Then we also assume that
$N(\bar{\rho})=N$. Finally $f$ is the normalized primitive form
corresponding to an elliptic curve $E$, which we assume has good
ordinary reduction at $p$.

As we have already mentioned in the introduction we will show that
our local ring $h_{2}(\Gamma_{0}(Npm^{2}),\ri)_{\m}$ is isomorphic
to some reduced ring $\T_{\Sigma}$ that appear in Wiles' work. For
the purposes of this section we do not really need to establish this
identification as we can work only with the full Hecke algebra.
However in order to make some remarks when we later consider the
case $p>3$, of course in a Hilbert modular form setting, we will
refer to this identification.

Let us start by fixing a general setting for this section. We write
$f$ for a normalized primitive form of conductor $N$ and trivial
character that is ordinary at $p$. We will write $\rho:=\rho_{f}$
for the $p$-adic representation that we attach to $f$, and
$\bar{\rho}:=\bar{\rho}_{f}$ for its mod $p$ representation. From
now on we assume that the local field $K$ is always sufficiently
large, in the sense that always contains all the Fourier
coefficients of the cusp forms that we consider. We fix a set
$\Sigma$ of primes $\ell \neq p$ that do not divide $N$. We define,
$N_{\Sigma}:=N p \prod_{\ell \in \Sigma}\ell^{2}$. Moreover we write
${\cal N}_{\Sigma}$ for the set of primitive forms $g$ of conductor
dividing $N_{\Sigma}$ and trivial character with the property that,
\[
a(q,g) \mod{\lambda'} = tr(\overline{\rho}(Frob_{q}))
\,\,\,\,\,\forall q\,\,\,\,with\,\,\,\,\,\,\, (q,N_{\Sigma})=1
\]
where $\lambda'$ is the maximal ideal in the field $K_{g}$, the
minimal $\Q_{p}$ extension that contains the coefficients of $g$.
Notice that if we write $\rho_{g}$ for the $p$-adic representation
that we can attach to the newform $g$ then the above condition gives
that $\overline{\rho} \otimes_{k} k_{g} \cong \overline{\rho}_{g}$
where $\overline{\rho}_{g}$ is the unique, up to isomorphism, mod
$p$ semi-simple representation such that
$tr(\overline{\rho}_{g}(Frob_{q}))= a(q,g) \mod \lambda'$ for all
$(q,N(\overline{\rho}_{g})p)=1$. Let $g \in {\cal N}_{\Sigma}$ be
any of our selected primitive forms and for any such $g$ we consider
the normalized eigenform $g' \in S_{2}(\Gamma_{0}(N_{\Sigma}))$ that
is defined by,

\begin{enumerate}
    \item $a(q,g')=a(q,g)$ if $q$ does not divide $N_{\Sigma}/N_{g}$
    \item $a(q,g')=0$ if $q \neq p$ and $q$ divides $N_{\Sigma}/N_{g}$
    \item $a(p,g'):=u(g)$, the $p$-adic unit root of the equation
    $X^{2}-a(p,g)X+p=0$ if $p$ divides $N_{\Sigma}/N_{g}$. Note that
    this makes sense as $\bar{\rho}$ comes from an elliptic curve with good ordinary reduction at $p$.
\end{enumerate}

It follows from lemma 4.6 of \cite{DDT} that the form $\overline{g'}
\in S_{2}(\Gamma_{0}(N_{\Sigma});k)$ is independent of $g$, and more
precisely it is characterized by the conditions (1)\,\,
$a(q,\overline{g'})=tr\overline{\rho}_{I_{q}}(Frob_{q})$ if $q=p$ or
$q \not \in \Sigma$, where we write $\overline{\rho}_{I_{q}}$ for
the coinvariant space of the inertia at $q$, and
(2)\,\,$a(q,\overline{g'})=0$ if $q \in \Sigma$. We then write $\m$
for the maximal ideal in $h_{2}(\Gamma_{0}(N_{\Sigma});\ri)$ that
corresponds to $\overline{g'}$ by proposition
\ref{proposition:structure}. Then we claim that
$h_{2}(\Gamma_{0}(N_{\Sigma});\ri)_{\m} \otimes_{\ri} K$ is a
semi-simple $K$-algebra. Indeed we have that,

\begin{proposition}
There is a $K$-algebra isomorphism,
\[
\phi : h_{2}(\Gamma_{0}(N_{\Sigma});\ri)_{\m} \otimes_{\ri} K \isom
\prod_{g \in \mathcal{N}_{\Sigma}} K
\]
given by
\[
(\phi(T_{q}))_{g} = \left\{%
\begin{array}{ll}
    a(q,g), & \hbox{if $q \not \in \Sigma \cup \{p\}$;} \\
    0, & \hbox{if $q \in \Sigma$;} \\
    u(g), & \hbox{if $q=p$.} \\
\end{array}%
\right.
\]
\end{proposition}

\paragraph{Proof}
By lemma \ref{lemma:Hecke algebras} we have an isomorphism
\[
h_{2}(\Gamma_{0}(N_{\Sigma});K) \isom \prod_{g \in
\mathcal{P}(N_{\Sigma})} A_{K,g}/I_{K,g}
\]
where $\mathcal{P}(N_{\Sigma})$ is the set of primitive forms in
$S_{2}(\Gamma_{0}(N_{\Sigma});K)$. Recall also that for a maximal
ideal $\m \subset  h_{2}(\Gamma_{0}(N_{\Sigma});\ri_{K})$ we have
that
\[
h_{2}(\Gamma_{0}(N_{\Sigma});\ri_{K})_{\m} \otimes_{\ri_{K}} K \isom
\prod_{\p} h_{2}(\Gamma_{0}(N_{\Sigma});K)_{\p}
\]
for all prime ideals $\p \subset h_{2}(\Gamma_{0}(N_{\Sigma});K)$
that restricted to $h_{2}(\Gamma_{0}(N_{\Sigma});\ri_{K})$ they are
contained in $\m$. Hence we obtain the following isomorphism,
\[
h_{2}(\Gamma_{0}(N_{\Sigma});\ri_{K})_{\m} \otimes_{\ri_{K}} K \isom
\prod_{g \in \mathcal{P}(N_{\Sigma})} \prod_{\p \in {\cal
M}_{g}}(A_{K,g}/I_{K,g})_{\p}
\]
where ${\cal M}_{g}$ is the set of prime ideals in $A_{K,g}/I_{K,g}$
that their image under the map $\Theta_{g}$ (with notation as in
lemma \ref{lemma:Hecke algebras}) when restricted to
$h_{2}(\Gamma_{0}(N_{\Sigma});\ri_{K})$ is in $\m$. But if $g$ is
not in ${\cal N}_{\Sigma}$ then ${\cal M}_{g}$ is empty. If $g$ is
in ${\cal N}_{\Sigma}$ then there is a unique prime ideal in ${\cal
M}_{g}$, call it $\p_{g'}$, that restricts inside $\m$. Indeed, it
is the prime ideal that corresponds to the normalized eigenform $g'$
constructed above as we have shown that it is the unique normalized
eigenform in $S_{K,g}$ with the required reduction. This prime ideal
is actually the kernel of the map $A_{K,g}/I_{K,g} \rightarrow K$
sending $u_{f,g} \mapsto a(q,g')$, and after localizing we obtain
$(A_{K,g}/I_{K,g})_{\p_{g'}} \isom K$. Finally the explicit
description of the isomorphism in lemma \ref{lemma:Hecke algebras}
gives the description of the isomorphism $\phi$. \qed

Now we are going to introduce the algebras $\T_{\Sigma}$ that appear
in Wiles' work. As we mentioned, we will not make any direct use of
them in the case of $p=3$. We consider the $\ri_{K}$-algebra
$\T'_{\Sigma} := \prod_{g \in {\cal N}_{\Sigma}} {\ri}_{K_{g}}$. and
we define the $\ri_{K}$-subalgebra $\T_{\Sigma} \subset
\T'_{\Sigma}$ generated over $\ri_{K}$ by the elements $T(q):=
(a(q,g))_{g}$ for all $q$ relatively prime to $N_{\Sigma}$. Note
that $\T_{\Sigma}$ is reduced as we consider only the ``good'', i.e.
away from the level, Hecke eigenvalues. In Wiles' work this algebra
is shown to be a deformation ring of the representation
$\bar{\rho}$. One of the crucial steps in his work is that he
identifies this algebra with a localized part of the full Hecke
algebra. In our case it follows from Proposition 4.7 in \cite{DDT}
that there is an isomorphism of $\ri_{K}$-algebras, $\phi:
h_{2}(\Gamma_{0}(N_{\Sigma});\ri_{K})_{\m} \isom \T_{\Sigma}$ given
by $T_{q} \mapsto T(q)$ for all $q$ relatively prime to
$N_{\Sigma}$.

\paragraph{Modules of congruences} For our primitive form $f \in S_{2}(\Gamma_{0}(N_{f};\Z_{p})$ we are
interested in the map
$\pi_{\Sigma}:h_{2}(\Gamma_{0}(N_{\Sigma});\Z_{p}) \rightarrow
\Z_{p}$ induced from the normalized eigenform $f' = \tilde{f}_{0} =
f_{0}|\imath_{m}$, with $m=\prod_{\ell \in \Sigma}\ell$. By the
universal property of localization, this map factors as,
\[
\pi_{\Sigma} : h_{2}(\Gamma_{0}(N_{\Sigma});\Z_{p}) \rightarrow
h_{2}(\Gamma_{0}(N_{\Sigma});\Z_{p})_{\m} \rightarrow \Z_{p}
\]
where the maximal ideal $\m$ is as in the previous section. If we
use the identification
$h_{2}(\Gamma_{0}(N_{\Sigma});\Z_{p})_{\m}\isom \T_{\Sigma}$ then
the map can be realized as the projection to the component
corresponding to $f$. Moreover we have shown that
$h_{2}(\Gamma_{0}(N_{\Sigma});\Z_{p})_{\m}\cong \T_{\Sigma}$ is
reduced and in particular the map $\pi_{\Sigma}$ induces a splitting
$h_{2}(\Gamma_{0}(N_{\Sigma});\Z_{p})_{\m} \otimes_{\Z} \Q_{p} =
\Q_{p} \times A$, where we write just $\Q_{p}$ as $f$ has rational
coefficients. Recall that we write $1_{\Q_{p}}$ for the idempotent
corresponding to the copy of $\Q_{p}$. We would like to study its
``denominator'' i.e. a quantity $c(f,m)$ such that $c(f,m)
1_{\Q_{p}}$ is integral. For this we now introduce the notion of the
module of congruences.

We start with some general definitions and properties of the module
of congruences as for example are given by Hida in his book, see
\cite{Hida9} (page 276). Let us write $h$ for a local ring
$h_{2}(\Gamma_{0}(N),\ri_{K})_{\m}$ for some $N$. Moreover let us
assume that $h$ is reduced and that we are given a map $\phi:h
\rightarrow \ri_{K}$, such that it induces an $K$-algebra
decomposition,
\[
h \otimes_{\ri_{K}} K \cong K \times A
\]
for some $K$-algebra $A$. Let us denote by $1_{\phi}$ the idempotent
that corresponds to the first summand $K$. We define
$\mathfrak{a}:=Ker(h\rightarrow A)$ and $\wp := Ker(\phi)$. Note
also that $Ann_{h}(\wp) = \mathfrak{a}$.

\begin{definition}
The module of congruences $C_{0}(h)$ of $\phi : h \rightarrow
\ri_{K}$ is defined as,
\[
C_{0}(h) := (h/\mathfrak{a}) \otimes_{h,\phi} \ri_{K} \cong
\frac{h}{\mathfrak{a} \oplus \wp} \cong \ri_{K}/\phi(\mathfrak{a})
\cong 1_{\phi}h/\mathfrak{a}
\]
\end{definition}

We now consider the module of congruences $C_{0}(h)$ for our reduced
ring $h:=h_{2}(\Gamma_{0}(N_{\Sigma}),\Z_{p})_{\m}$ and our map
$\pi_{\Sigma}$. We will compare it with a ``cohomological'' module
of congruences, following the terminology of Hida and Ribet, which
we will introduce below. Let us write $X$ for the compact modular
curve $X_{0}(N_{\Sigma})$. We consider the first cohomology group
$H^{1}(X,\Z_{p})$ and we have seen in chapter three that this as a
Hecke module over $h_{2}(\Gamma_{0}(N_{\Sigma});\Z_{p})$. Moreover
we consider the standard skew-symmetric bilinear perfect pairing as
for example in \cite{DDT} page 106,
\[
(\cdot\,,\cdot) : H^{1}(X;\Z_{p}) \times H^{1}(X;\Z_{p}) \rightarrow
\Z_{p}
\]
Let us define $L:=H^{1}(X;\Z_{p})_{\m}$, where we have localized
$H^{1}(X;\Z_{p})$ at the maximal ideal $\m$. This is then an
$h$-module. Let us consider the action of complex conjugation on
$H^{1}(X;\Z_{p})$ and define $L[+],L[-]$ for the eigenspaces of $L$.
We define a cohomological module of congruences by,
\[
C^{coh}(L[+]):=1_{\Q_{p}}L[+]/1_{\Q_{p}}L[+] \cap L[+] \cong
L[+]^{\Q_{p}}/L[+]_{\Q_{p}}\cong \frac{L[+]}{L[+][\wp] \oplus
L[+][\mathfrak{a}]}
\]
where we have set $L[+]^{\Q_{p}}:=1_{\Q_{p}}L[+]$, the projection of
$L[+]$ to the first component of the decomposition
$L[+]\otimes_{\Z_{p}} \Q_{p} =1_{\Q_{p}}(L[+] \otimes_{\Z_{p}}
\Q_{p}) \oplus (1-1_{\Q_{p}})(L[+] \otimes_{\Z_{p}} \Q_{p})$ induced
by the splitting of the Hecke algebra. Also
$L[+]_{\Q_{p}}:=1_{\Q_{p}}L[+] \cap L[+] = L[+][\wp]$ with $\wp$ and
$\mathfrak{a}$ as in the definition of the module of congruences,
the restriction of $L[+]$ to the first component. Here, as usual,
$L[+][\wp]=\{ l \in L[+] : \lambda(l) = 0,\,\,\, \forall \lambda \in
\wp\}$. Note that $L \otimes_{\Z_{p}} \Q_{p} = H^{1}(X;\Z_{p})_{\m}
\otimes_{\Z_{p}} \Q_{p} $ is free of rank two over $h
\otimes_{\Z_{p}} \Q_{p}$ and hence $L[\pm] \otimes_{\Z_{p}} \Q_{p}$
of rank one. In particular we have that $L[+]_{\Q_{p}}$ is a free
$\Z_{p}$-module of rank one. We fix a basis $x_{+}$. The same holds
for $L[-]_{\Q_{p}}$ and we fix a basis $x_{-}$.

\begin{lemma} For the ``cohomological'' module of congruences we
have,
\[
C^{coh}(L[+]) \cong \Z_{p}/((x_{+},x_{-}))
\]
\end{lemma}

\paragraph{Proof} (See also \cite{DDT} p. 105 and \cite{Hida9} p. 275) It is enough to show
that $L[+]_{\Q_{p}} \cong Hom_{\Z_{p}}(L[-]^{\Q_{p}},\Z_{p})$ and
$L[+]^{\Q_{p}} \cong Hom_{\Z_{p}}(L[-]_{\Q_{p}},\Z_{p})$. This
follows by the prefect pairing $(\cdot,\cdot)$ on $L$. Indeed first
we note that complex conjugation acts as
$(a,b^{\rho})=-(a^{\rho},b)$, which explains the eigenspaces. Now
let us show that $L[+]^{\Q_{p}} \cong
Hom_{\Z_{p}}(L[-]_{\Q_{p}},\Z_{p})$ as the other claim is obtained
similarly. Note that if we consider a basis $\{x_{1},x_{2},\cdots\}$
of $L$ as a $\Z_{p}$ module such that $x_{1}=x_{-}$ then as $L$ is
self-dual with respect to $(\cdot,\cdot)$ there is a dual basis
$\{x_{-}^{*},x_{2}^{*},\cdots \}$ in $L$. Taking the projection
$1_{\Q_{p}}x_{-}^{*}[+]$ gives a dual basis of $L[-]_{\Q_{p}}$. \qed

We would like to compare the module of congruences $C_{0}(h)$ and
$C^{coh}(L[+])$. Under our assumptions on $f$, i.e. it is
$p$-ordinary and its modulo $p$ representation is irreducible, we
have the following important theorem of Wiles \cite{Wiles}.

\begin{theorem}$($Wiles$)$ The $h$-module $H^{1}(X,\Z_{p})_{\m}$ is
free (of rank two).
\end{theorem}
We can conclude,
\begin{corollary} $C_{0}(h) \cong \frac{h}{\wp \oplus \mathfrak{a}} \cong \frac{1_{\Q_{p}}h}{\mathfrak{a}} \cong L[+]^{\Q_{p}}/L[+]_{\Q_{p}} \cong
C^{coh}(L[+]) \cong \Z_{p}/((x_{+},x_{-}))$.
\end{corollary}

Hence the quantity $(x_{+},x_{-}) \in \Z_{p}$ annihilates the module
of congruences and in particular we know that
$(x_{+},x_{-})1_{\Q_{p}} \in h \subset
h_{2}(\Gamma_{0}(N_{\Sigma}),\Z_{p})$. Hence we can define, up to
$p$-adic units, $c(f,m):=(x_{+},x_{-})$. In the next section we will
study the relation of $c(f,m)$ with the periods
$<\tilde{f}_{0}|\tau_{N_{\Sigma}},\tilde{f}_{0}>$ that appear in the
first form of our congruences and eventually relate it to the
N\'{e}ron periods $\Omega_{+}(E)$ and $\Omega_{-}(E)$.

\paragraph{Relations between periods of different levels:}The main aim
now is to understand the relation between the quantity $c(f,m)$ and
the automorphic periods
$<\tilde{f}_{0}|\tau_{N_{\Sigma}},\tilde{f}_{0}>$ that appear in our
congruences in theorem \ref{theorem:Congruences I}. Recall that we
are considering the homomorphism $\pi_{\Sigma}:
h_{2}(N_{\Sigma};\Z_{p})_{\m} \rightarrow \Z_{p}$ corresponding to
our normalized eigenform $g:=\tilde{f}_{0}$, arising from the
newform $f$ (i.e. $g=f_{0}|\imath_{m})$. Let us write $\Z_{(p)}$ for
the localization of $\Z$ at $p$. Recall that we write $\wp$ for the
kernel of $\pi_{\Sigma}$. Then we have the inclusion,
\[
H^{1}(X_{0}(N_{\Sigma});\Z_{(p)})[g] \subset
H^{1}(X_{0}(N_{\Sigma});\Z_{p})_{\m}[\wp]=L[\wp]
\]
Let us choose a basis $\{x_{+},x_{-}\}$ for $L[\wp]$ which is in the
image of $H^{1}(X_{0}(N_{\Sigma});\Z_{(p)})[g]$. We remind the
reader that $p=3$ and so as we are interested in statements up to
$p$-adic units we can keep working with eigenspaces. We consider the
$\C$ vector space $H^{1}(X_{0}(N_{\Sigma});\C)[g]$. The classical
Eichler-Shimura isomorphism gives,
\[
S_{2}(\Gamma_{0};\C) \oplus \overline{S_{2}}(\Gamma_{0};\C)
\stackrel{\sim}{\rightarrow} H^{1}(X_{0}(N_{\Sigma});\C)
\]
where we write $\overline{S_{2}}(\Gamma_{0}(N_{\Sigma});\C)$ for the
space of the anti-holomorphic cusp forms. A canonical basis of
$H^{1}(X_{0}(N_{\Sigma});\C)[g]$ is given by,
$\{\omega_{g},\overline{\omega_{g^{\rho}}}\}$ where $\omega_{g}=\sum
a(n,g)q^{n-1}dq$ is a holomorphic differential on $X$  and
$\overline{\omega}_{g^{\rho}}=\sum a(n,g)\bar{q}^{n-1}d\bar{q}$ an
anti-holomorphic. Note that actually in the case of interest $g$ has
rational coefficients, hence $g^{\rho}=g$. We now define a period
$\Omega(f)_{\Sigma}$ as follows. We let $A_{\Sigma}$ be the two by
two invertible matrix in $GL_{2}(\C)$ such that
$[\overline{\omega}_{g},\omega_{g}]=[x_{+},x_{-}]A_{\Sigma}$ and we
define $\Omega(f)_{\Sigma}:=det(A_{\Sigma})$. Then,

\begin{lemma}\label{lemma:relation} With notation as above we have the following
equation,
\[
c(f,m) \Omega(f)_{\Sigma} = (x_{+},x_{-}) \Omega(f)_{\Sigma} =
<g|\tau_{N_{\Sigma}},g>
\]
\end{lemma}

\paragraph{Proof}Just note that the skew-symmetry of the pairing for $A_{\Sigma} = \left(%
\begin{array}{cc}
  a & b \\
  c & d \\
\end{array}%
\right)$ gives $(\overline{\omega}_{g},\omega_{g}) = (ax_{+} +
cx_{-},bx_{+}+dx_{-}) =
ad(x_{+},x_{-})-cb(x_{+},x_{-})=det(A_{\Sigma})(x_{+},x_{-})$. But
by the definition of the pairing we have that $
(\overline{\omega}_{g},\omega_{g}) =
\int_{X}\overline{\omega}_{g|\tau_{N_{\Sigma}}}\wedge \omega_{g} =
<g|\tau_{N_{\Sigma}},g>$  \qed

Hence in view of theorem \ref{theorem:Congruences I}, in order to
conclude the congruences we need to relate the automorphic periods
$\Omega(f)_{\Sigma}$ with the periods $\Omega_{+}(E)\Omega_{-}(E)$.
The following theorem is taken from \cite{DDT}, p.108 and is based
on Wiles' generalization of the so-called Ihara's lemma.

\begin{theorem}\label{theorem:periodsLevel} We have, up to $p$-adic units,
$\Omega(f)_{\Sigma}=\Omega(f)_{\Sigma=\emptyset}=\Omega(f)$ where
$\Omega(f)$ is the period defined by taking $g=f$ above.
\end{theorem}

Now we relate the period $\Omega(f)$ with the N\'{e}ron periods
$\Omega_{+}(E)\Omega_{-}(E)$, see \cite{Hida1}, p. 255 and
\cite{Wiles} p.537. Recall that we consider an elliptic curve $E/\Q$
of conductor $N$ and we write $f \in S_{2}(\Gamma_{0}(N);\Z)$ for
the primitive form of weight two and conductor $N$ associated to it.
Let us fix a global minimal Weierstrass equation of $E$ over $\Z$
and let denote by $\omega_{E}$ the N\'{e}ron differential of this
equation. Let us consider the eigenspaces of $H_{1}(E(\C),\Z)$ under
the action of complex conjugation. We fix generators
$\gamma_{E}^{+}$ and $\gamma_{E}^{-}$ for the spaces
$H_{1}(E(\C);\Z)^{+}$ and $H_{1}(E(\C);\Z)^{-}$. Recall that we have
defined the N\'{e}ron periods as,
\[
\Omega(E)_{\pm} := \int_{\gamma_{E}^{\pm}} \omega_{E}
\]
If we write $\phi:X_{0}(N)\rightarrow E$ for the strong Weil
parametrization of $E/\Q$ then we have that $\phi^{*}\omega_{E} = 2
\pi \imath c_{E} f(z)dz$ where $c_{E} \in \Q^{\times}$ and in
particular it has been proved by Mazur \cite{Mazur2} that $c_{E}$ is
a $p$-adic unit if $p^{2} \nmid 4N$. By Poincar\'{e} duality we can
pick a $\Z_{(p)}$ basis $\{c_{1}, c_{2}, \ldots, c_{m}\}$ of
$H_{1}(X_{0}(N);\Z_{(p)})$ such that $\int_{c_{j}}\ell_{i} =
\delta_{ij}$ for $i,j=\{1,2\}$ and $\ell_{1},\ell_{2}$ a basis of
$H^{1}(X_{0}(N);\Z_{(p)})[\p]$. Then we have,

\[
det\left(%
\begin{array}{cc}
  \int_{c_{1}}\omega_{f} & \int_{c_{1}}\bar{\omega}_{f} \\
  \int_{c_{2}}\omega_{f} & \int_{c_{2}}\bar{\omega}_{f} \\
\end{array}%
\right) =  det\left(%
\begin{array}{cc}
  \int_{c_{1}}\ell_{1} & \int_{c_{1}}\ell_{2} \\
  \int_{c_{2}}\ell_{1} & \int_{c_{2}}\ell_{2} \\
\end{array}%
\right) \Omega(f) = \Omega(f)
\]
Moreover we have,
\[
\left| det \left(%
\begin{array}{cc}
  \int_{\phi(c_{1})}\omega & \int_{\phi(c_{1})}\bar{\omega} \\
  \int_{\phi(c_{2})}\omega & \int_{\phi(c_{2})}\bar{\omega} \\
\end{array}%
\right) \right| = 4 \pi^{2} c_{E}^{2} \left| det \left(%
\begin{array}{cc}
  \int_{c_{1}}\omega_{f} & \int_{c_{1}} \bar{\omega}_{f}\\
  \int_{c_{2}}\omega_{f} & \int_{c_{2}} \bar{\omega}_{f} \\
\end{array}%
\right) \right|=4 \pi^{2} c_{E}^{2} \Omega(f)
\]
And also,
\[
\left| det\left(%
\begin{array}{cc}
  \int_{\gamma_{E}^{+}}\omega & \int_{\gamma_{E}^{+}}\bar{\omega} \\
  \int_{\gamma_{E}^{-}}\omega & \int_{\gamma_{E}^{-}}\bar{\omega} \\
\end{array}%
\right)\right| = 2|\Omega_{+}(E)\Omega_{-}(E)| =
2\Omega_{+}(E)\Omega_{-}(E)\imath^{-1}
\]
where the last equality follows from the fact that we can pick
$\gamma_{E}^{\pm}$ such that $\imath^{-1} \Omega_{-}(E)$ and
$\Omega_{+}(E)$ are real positive. As
\[
\left| det\left(%
\begin{array}{cc}
  \int_{\gamma_{E}^{+}}\omega & \int_{\gamma_{E}^{+}}\bar{\omega} \\
  \int_{\gamma_{E}^{-}}\omega & \int_{\gamma_{E}^{-}}\bar{\omega} \\
\end{array}%
\right)\right| = \left| det \left(%
\begin{array}{cc}
  \int_{\phi(c_{1})}\omega & \int_{\phi(c_{1})}\bar{\omega} \\
  \int_{\phi(c_{2})}\omega & \int_{\phi(c_{2})}\bar{\omega} \\
\end{array}%
\right) \right|
\]
up to $p$-adic units, we have,

\begin{theorem}The relation of the period $\Omega(f)$ with
the periods $\Omega_{+}(E)$ and $\Omega_{-}(E)$, up to $p$-adic
units, is given by the equation,
\[
\pi^{2} \imath \Omega(f) = \Omega_{+}(E)\Omega_{-}(E)
\]
\end{theorem}

Putting all together, theorem \ref{theorem:Congruences I}, lemma
\ref{lemma:relation}, theorem \ref{theorem:periodsLevel} and the
above theorem we conclude

\begin{theorem}
Consider an elliptic curve $E$ as in the introduction. Let $m$ be a
power free positive integer with $(m,N_{E})=(m,p)=1$ with $p=3$.
Consider the Galois extension $\Q(\mu_{p},\sqrt[p]{m})/\Q$ and let
$\rho$ be the unique non-trivial two dimensional
Artin-representation that factors through
$Gal(\Q(\mu_{p},\sqrt[p]{m})/\Q)$. Then,

\[
e_{p}(\rho)u^{-v_{p}(N_{\rho})}\frac{P_{p}(\rho,u^{-1})}{P_{p}(\rho,w^{-1})}\frac{L_{\{p,q|m\}}(E\otimes
\rho,1)}{\Omega_{+}(E)\Omega_{-}(E)} \equiv
e_{p}(\sigma)u^{-v_{p}(N_{\sigma})}\frac{P_{p}(\sigma,u^{-1})}{P_{p}(\sigma,w^{-1})}\frac{L_{\{p,q|m\}}(E\otimes
\sigma,1)}{\Omega_{+}(E)\Omega_{-}(E)} \mod{p}
\]

where $\sigma = 1 \oplus \epsilon_{p}$ with $\epsilon_{p}$ the
non-trivial character of $\Q(\mu_{p})/\Q$ and $u,w$ such that,
\[
1-a_{p}X+pX^{2} = (1-uX)(1-wX),\,\,\,\, u\in
\Z_{p}^{\times}\,\,\,and\,\,\,p+1-a_{p}=\#E_{p}(\mathbb{F}_{p})
\]
\end{theorem}

Let us remark here that it is easy to see that we could relax our
assumption that the conductor of $E$ equals the conductor of the mod
$p$ representation in the expense of obtaining the weaker
congruences,
\[
R(\rho)\prod_{q|N_{diff}}P_{q}(E,\rho,1) \equiv R(\sigma)
\prod_{q|N_{diff}}P_{q}(E,\sigma,1) \mod p
\]
where $N_{diff}$ the ratio of the conductor of $E$ over the Artin
conductor of the mod $p$ representation. Indeed, instead of
considering the eigenform $f_{0}|\imath_{m}$ we need to consider the
one where we remove the primes that divide $m$ and those that divide
$N_{diff}$, i.e. $f_{0}|\imath_{mN_{diff}}$. It is this eigenform
that will induce a homomorphism of the Hecke algebra that factors
through a reduced local ring in case that $N_{diff}$ is not one.
Then everything carries as above but eventually we remove also the
Euler factors at $N_{diff}$ as we have modified $f$ in this way.

\section{Speculations for the case $p>3$}

As the title indicates there are no real results in this section.
The aim is to give a brief account of the problems that we face
trying to extend our previous results to the case $p>3$, working in
the Hilbert modular form setting.

Note that in the previous section the fundamental result of Wiles
allowed us to compare the size of the module of congruences for the
Hecke algebra with the cohomological one which eventually was
related to the periods that we used. In particular the crucial
results were, first, that the localized first cohomology group was a
free Hecke module over the local ring corresponding to our cusp form
and second the generalization of ``Ihara's lemma'' that allowed us
to relate the different levels. In the Hilbert modular form setting,
results of this form have be obtained by Diamond in \cite{Diamond1}
and Dimitrov \cite{Dimitrov1}. However for both authors it is
crucial to assume that they work with a prime that is unramified in
the totally real field.

Moreover there is another difficulty that is related to the
automorphic periods that we can also define in this Hilbert modular
forms setting. Indeed using the Eichler-Shimura-Harder isomorphism
we can define periods $\Omega(\phi)_{\Sigma}$, the analogue of
$\Omega(f)_{\Sigma}$, and their relation to the Petersson inner
product is governed by the cohomological module of congruences.
However even if we had an Ihara type lemma in this case we would
have still to relate $\Omega(\phi)$, the minimal level, to the
N\'{e}ron periods up to $p$-adic units. So we run again into the
same question as the one we addressed in our work \cite{Boug}, that
is to understand the behavior of the automorphic periods under base
change.

Having stated these problems we would like to speculate a little.
Note that in what we said above we do not really make use of the
fact that actually we consider a Hilbert modular form that is coming
from base-change. In what follows we will try to indicate that
perhaps one can avoid working over the totally real field and reduce
our questions to the study of the adjoint square $L(ad(f),s)$
$L$-function associated to $f$ and its behavior under twists over
the extension $F/\Q$, $F=\Q(\mu_{p})^{+}$. This also will justify
our choice to underline the identification of the local ring
$h_{2}(\Gamma_{0}(N_{\Sigma}),\Z_{p})_{\m}$ with the ``deformation''
ring $T_{\Sigma}$, in the previous section. Our exposition is very
brief and not rigorous.

So we keep the same notation as in the previous sections with the
obvious extensions to the Hilbert modular case. That is $h$ is now a
local ring of the Hecke algebra acting on the space of Hilbert cusp
forms of level $N\p m^{2}$ completed at $p$. Moreover it is the
reduced local ring through which our ordinary normalized cusp form $
\tilde{\phi}_{0}$ factors. We write $C_{0}(h)$ for its module of
congruences. As in the elliptic case one can identify $h$ with the
``deformation'' ring $\T_{\Sigma}$. We consider the
$\Z_{(p)}$-module, $C_{1}(h):=(ker
\pi_{\Sigma})/(ker\pi_{\Sigma})^{2}$. Then by \cite{DDT} page 117,
we have the inequality $|C_{1}(h)| \geq |C_{0}(h)|$. Let us now
write $\rho_{F}$ for the 3-adic representation obtained from $\rho$
by restriction to $G_{F}$ and consider its reduction
$\bar{\rho}_{F}$  modulo 3. Now let us impose the following
conditions on $\rho_{F}$,

\begin{enumerate}

\item $\bar{\rho}_{F}$ is absolutely irreducible,

\item ($p$-ordinary) $\rho_{F}|_{D_{\p}} \cong \left(%
\begin{array}{cc}
  \delta_{\p} & * \\
  0 & \epsilon_{\p} \\
\end{array}%
\right)$ where $\delta_{\p}$ is an unramified character.
\end{enumerate}

Then it is known by the work of Mazur that there exists a universal
deformation couple $(R_{F},\varrho_{F})$, in the terminology of Hida
\cite{Hida9}, that represents deformations with prescribed
determinant and ramification in a way that we do not make explicit
here. As we indicated in the elliptic case, the algebra
$\T_{\Sigma}$ can be interpreted as a deformation algebra for
$\bar{\rho}_{F}$ and hence there is a surjection $R_{F}
\twoheadrightarrow \T_{\Sigma}$. This implies the inequality
\cite{DDT} p.118, $|C_{1}(R_{F})| \geq |C_{1}(\T_{\Sigma})| \geq
|C_{0}(\T_{\Sigma})|$. Again by Mazur's theory one can identify
$C_{1}(R_{F})$ with the Pontryagin dual of a properly defined Selmer
group $Sel(ad(\rho_{F}))$ attached to $ad (\rho_{F})$. But we can
decompose $Sel(ad(\rho_{F})) = \oplus_{\chi}Sel(ad(\rho)\otimes
\chi)$ for $\chi \in Gal(F/\Q)^{V}$. So back to our congruences we
have a bound for our constant $c(\phi,m)$ by the the sizes of the
Selmer groups of $ad(\rho)$ twisted by characters $\chi$ that factor
through $Gal(F/\Q)$.

Recall that the periods that appear in our congruences involve the
Petersson inner product
$<\tilde{\phi}_{0}|\tau_{N_{\Sigma}},\tilde{\phi}_{0}>$. We would
like to factor this quantity to quantities that are related with the
cusp form $f$ and more important we would like to obtain some
control of the constants that may appear. In the elliptic modular
forms case a formula of Shimura allows one to relate the Petersson
inner product $<\tilde{f}_{0}|\tau_{N_{\Sigma}},\tilde{f}_{0}>$ to
the value of the adjoint $L$ function $L(ad(f),s)$ at $s=2$, in
particular they are equal up to powers of $\pi$ and modified Euler
factors at primes dividing $N_{\Sigma}$. This formula can be
extended to the Hilbert modular case, see \cite{Shimura4} page 669.
One can then use the inductive properties of the $L$ functions to
rewrite the periods
$<\tilde{\phi}_{0}|\tau_{N_{\Sigma}},\tilde{\phi}_{0}>$ as a product
of the form $\prod_{\chi}L(ad(f)\otimes \chi,2)$, up to modified
Euler factors and powers of $\pi$. These modified Euler factors
should correspond to the local conditions that we have impose to the
above mentioned Selmer group depending on the deformation problem.

Granted all the above speculations, we see that proving the
congruences for $p>3$ is closely related to the Tamagawa number
conjecture for the adjoint square $L$ function of $f$ and its twists
with characters over the extension $Gal(F/\Q)$.

\end{document}